\documentclass[11pt]{amsart}
\usepackage[matrix,arrow,curve]{xy}
\usepackage{amssymb}
\calclayout
\newdir{  >}{{}*!/-10pt/\dir{>}}
\newtheorem{dummy}{realdumb}[section]
\newtheorem{theorem}[dummy]{Theorem}
\newtheorem*{thma}{Theorem A}
\newtheorem*{thmb}{Theorem B}
\newtheorem*{thmc}{Theorem C}
\newtheorem*{thmd}{Theorem D}
\newtheorem*{thme}{Theorem E}
\newtheorem*{thmmain}{Theorem}
\newtheorem{lemma}[dummy]{Lemma}
\newtheorem{corollary}[dummy]{Corollary}
\newtheorem{proposition}[dummy]{Proposition}
\theoremstyle{definition}               
\newtheorem{definition}[dummy]{Definition}
\newtheorem{example}[dummy]{Example}
\newtheorem{remark}[dummy]{Remark}
\newtheorem*{acknowledgement}{Acknowledgement}
\DeclareMathOperator{\Res}{Res}
\DeclareMathOperator{\Mod}{mod}
\DeclareMathOperator{\Ind}{Ind}
\DeclareMathOperator{\Inf}{inf}
\DeclareMathOperator{\Fix}{Fix}
\DeclareMathOperator{\Top}{Top}
\DeclareMathOperator{\wh}{Wh}
\DeclareMathOperator{\Iso}{Iso}
\DeclareMathOperator{\Image}{Im}
\newcommand{\trf}[1]{tr{\hskip -1.4truept}f_{#1}}
\newcommand{\cy}[1]{\bz/{#1}}
\newcommand{\Cy}[1]{C(#1)}
\newcommand{\adic}[2]{\widehat{#1}_{#2}}
\newcommand{\Zadictwo}[1]{\adic{\bz}{2}{#1}}
\newcommand{\zadictwo}{\adic{\bz}{2}}
\newcommand{\zlocal}{\widehat{\bz}}
\newcommand{\qlocal}{\widehat\bq}
\newcommand{\Lgroup}[3]{L^{#1}_{#2}({#3})}
\newcommand{\Godd}{G_{odd}}
\newcommand{\Kodd}{{F}}
\newcommand{\Wmax}{W_{max}}
\newcommand{\rtop}[1]{R_{\Top}(#1)}
\newcommand{\Rtop}{R_{\Top}}
\newcommand{\Rtopfree}{R^{\scriptscriptstyle free}_{\Top}}
\newcommand{\rtopfree}[1]{\Rfree_{\Top}(#1)}
\newcommand{\rtildetopfree}[1]{\Rtildefree_{\Top}(#1)}
\newcommand{\Rtildetopfree}{\Rtildefree_{\Top}}
\newcommand{\rfree}[1]{\Rfree(#1)}
\newcommand{\Rfree}{R^{\scriptscriptstyle free}}
\newcommand{\rtildefree}[1]{\Rtildefree(#1)}
\newcommand{\Rtildefree}{{\widetilde R}^{{\scriptscriptstyle free}}}
\newcommand{\orient}{\{\pm 1\}}
\newcommand{\rminus}{{\br_{-}}}
\newcommand{\Rminus}[1]{{\br^{#1}_{-}}}
\newcommand{\Rplus}[1]{{\br^{#1}_{+}}}
\newcommand{\rplus}{{\br_{+}}}
\newcommand{\Linfty}{L^{\scriptscriptstyle \langle -\infty \rangle}}
\newcommand{\Lminus}[1]{L^{\scriptscriptstyle \langle -{#1} \rangle}}
\newcommand{\Kminus}[1]{K_{\scriptscriptstyle -{#1}}}
\newcommand{\Ccat}[3]{\cc_{#1}^{>#2} ({#3})}
\newcommand{\ccat}[2]{\cc_{#1}({#2})}
\newcommand{\bcat}[1]{\cc_{#1,G}(\mathbf Z )}
\newcommand{\bcatrminus}[1]{\bcat{#1\times\rminus}}
\newcommand{\bCatrminus}{\bcat{\rminus}}
\newcommand{\bCat}[3]{\cc^{>#3}_{#1,#2}\,}
\newcommand{\brelcat}[2]{\Ccat{#1, G}{#2}{\mathbf Z}}
\newcommand{\brelcatrminus}[2]{\brelcat{#1\times\rminus}{#2\times\rminus}}
\newcommand{\brelCatrminus}[1]{\brelcat{#1\times \rminus}{\rminus}}
\newcommand{\Uell}{\bigcup \ell_\alpha}
\newcommand{\Vertical}{\,|\,}
\newcommand{\trfrminus}{\trf{\rminus}}
\newcommand{\ck}{{\mathbf k}}
\newcommand{\bz}{\mathbf Z}
\newcommand{\bc}{{\mathbf C}}
\newcommand{\bu}{\scriptscriptstyle\mathbf U}
\newcommand{\bU}{\mathbf U}
\newcommand{\bq}{\mathbf Q}
\newcommand{\br}{\mathbf R}
\newcommand{\bY}{\scriptscriptstyle \mathbf Y}
\newcommand{\cmlocal}{\widehat{\mathcal M}}
\newcommand{\BL}{{\mathbb L}}
\newcommand{\BK}{{\mathbb K}}
\newcommand{\cm}{{\mathcal M}}
\newcommand{\ca}{{\mathcal A}}
\newcommand{\bd}{\partial}
\newcommand{\cc}{{\mathcal C}}
\newcommand{\cu}{{\mathcal U}}
\newcommand{\ZG}{{\mathbf Z} G}
\newcommand{\ZK}{{\mathbf Z} K}
\newcommand{\ZH}{\mathbf Z H}
\newcommand{\QG}{\mathbf Q G}
\newcommand{\RA}[1]{\xrightarrow{#1}}
\newcommand{\rel}[2]{{#1}\hskip-2.2truept \to \hskip-2.2truept{#2}}

\newcommand{\Zreltwo}[1]{{\mathbf Z}{#1}\hskip-2.2truept \to \hskip-2.2truept\adic{\mathbf Z}{2}{#1}}
\newcommand{\Coh}[2]{H^{#1}(#2)}
\newcommand{\Units}[1]{{#1}^{\scriptscriptstyle\times}}
\newcommand{\localunits}[3]{\hat{#1}^{\scriptscriptstyle{\times  #2}}_{#3}}

\newcommand{\eqncount}{\setcounter{equation}{\value{dummy}}%
\addtocounter{dummy}{1}}
\newcommand{\partone}{\textbf{I}}
\newcommand{\parttwo}{\textbf{II}}
\newcommand{\pone}{[\textbf{I}]}
\newcommand{\ptwo}{[\textbf{II}]}

\begin{document}
  \title[Similarities of cyclic groups]{Topological equivalence of linear representations\\
 for cyclic groups: \parttwo}
\author{Ian Hambleton}
\address{Department of Mathematics \& Statistics,
McMaster University\newline
\indent Hamilton, Ont.,  Canada,  L8S 4K1}
\email{ ian@math.mcmaster.ca}

\author{Erik K. Pedersen}
\address{Department\ of\ Mathematical Sciences\\
SUNY at Binghamton\newline
\indent Binghamton,  NY, 13901}
\email{ erik@math.binghamton.edu}
\thanks{Partially supported by NSERC grant A4000
and  NSF grant  DMS 9104026. The authors also
wish to thank the Max Planck Institut f\"ur Mathematik, Bonn, for
its hospitality and support.}
\begin{abstract} In  the two parts of this paper we prove that the Reidemeister torsion invariants
determine  topological equivalence of $G${--}representations,
for $G$ a finite cyclic group.
\end{abstract}
\date{Feb. 19. 2004}
\maketitle
\section{ Introduction\label{intro}}
Let $G$ be a finite group and $V$, $V'$ finite dimensional
real orthogonal representations of $G$.
Then $V $ is said to be
\emph{topologically equivalent} to $V'$
(denoted $V \sim_t V'$) if there exists a homeomorphism $h\colon
V \to V'$ which is $G$-equivariant.
If $V$, $V'$ are topologically equivalent, but not linearly isomorphic,
then such a homeomorphism is called a non-linear similarity.
These notions were introduced and studied by de Rham \cite{drh2},
\cite{drh1}, and developed extensively in
 \cite{cs1}, \cite{cs2}, \cite{hsp1},  \cite{mr1}, and
\cite{cssw1}.
In the two parts of this paper, referred to as \pone\
and \ptwo, we complete de Rham's program by showing that
Reidemeister torsion invariants and number theory determine
non-linear similarity for finite cyclic groups.

A  $G$-representation is called \emph{free} if
each element $1\neq g \in G$ fixes only the zero vector.
Every representation of a finite cyclic group has a
unique maximal free subrepresentation.

\begin{thmmain} Let $G$ be a finite cyclic group
and $V_1$, $V_2$ be free
$G${--}representations. For any $G${--}representation $W$, the
existence of a non-linear similarity
$V_1\oplus W \sim_t V_2\oplus W$ is entirely 
determined by explicit
congruences in the weights of the free summands $V_1$, $V_2$,
 and the ratio
$\Delta(V_1)/\Delta(V_2)$ of their Reidemeister torsions, 
up to an algebraically described indeterminacy.
\end{thmmain}
\noindent
The notation and the  indeterminacy are given in Section
\ref{results} and a detailed statement of results in
Theorems A{--}E. 
This part of the paper contains the foundational results
and calculations in 
bounded algebraic $K$- and $L$-theory needed to prove
the main results on non-linear similarity. 
The study of non-linear similarities $V_1\oplus W \sim_t V_2\oplus W$ increases in difficulty with the number of isotropy types in
$W$. We introduce a new method using excision in bounded surgery theory, based on the \emph{orbit type filtration}, to organize and deal with these difficulties.
We expect that this technique will be useful for other applications.
Our most general results about non-linear similarity for arbitrary cyclic groups are Theorem C and its extensions (see
 Sections \ref{C} and \ref{rminustwo}).

In Sections \ref{splitting} and \ref{sD} we study the group
$\rtop{G}$ of $G${--}representations modulo \emph{stable} topological equivalences  (see \cite{cs2} where $\rtop{G}\otimes \bq$ is computed).
As an application of our general results,  we determine the
structure of the torsion in $\rtop{G}$, for $G$ any cyclic group
(see Theorem \ref{oneD}), and in Theorem D we give the  calculation of $\rtop{G}$
for $G = \Cy{4q}$, for $q$ odd, correcting \cite[Thm. 2]{cssw1}.
One interesting feature is that Corollary \ref{sixdim}
and  Theorem D indicate a connection between
 the orders of the \emph{ideal class groups} for
cyclotomic fields and topological equivalence of linear
representations.  

\tableofcontents
\section{ Statement of Results\label{results}}
For the reader's convenience, we  recall some notation from Part \partone, and then give the main results of both parts.
Theorems A and B  are proved in Part \partone\ and Theorems
C and D are proved in Part \parttwo. The proof of Theorem E
is divided between the two parts.

Let $G = \Cy{4q}$, where $q >1$, and let $H = \Cy{2q}$ denote
the subgroup of index 2 in $G$.
The maximal odd order subgroup of $G$ is denoted $\Godd$.
We fix a generator $G=\langle t \rangle$ and a primitive $4q^{th}$-root of
unity $\zeta = \exp{2\pi i/{4q}}$.
The group $G$ has both a trivial $1$-dimensional real representation,
denoted $\rplus$, and a non-trivial $1$-dimensional
real representation, denoted $\rminus$.

A \emph{free} $G${--}representation is a sum of faithful
1-dimensional complex representations.
Let $t^a$, $a \in \mathbf Z$, denote the complex numbers
 $\bc$ with action
$t\cdot z = \zeta^az$ for all $z \in \bc$.
This representation is free if and only if $(a, 4q) = 1$, and
 the coefficient $a$ is well{--}defined only modulo $4q$.
Since $t^a \cong t^{-a}$ as real $G${--}representations,  we
can  always choose the weights $a \equiv 1 \Mod{4}$.
This will be assumed unless otherwise mentioned.

Now suppose that $V_1 = t^{a_1} + \dots + t^{a_k}$
is a free $G${--}representation.
The Reidemeister torsion invariant of $V_1$ is defined
as
\[\Delta(V_1) = \prod_{i=1}^{k} (t^{a_i} -1) \in \mathbf Z [t]/\{\pm t^m\} \ .\]
Let $V_2 = t^{b_1} + \dots + t^{b_k}$ be another free representation,
such that $S(V_1)$ and $S(V_2)$ are $G$-homotopy equivalent. 
This just means that the products of the weights $\prod a_i \equiv \prod b_i
\Mod{4q}$.
Then the Whitehead torsion of any $G$-homotopy equivalence  is
determined by the element
\[\Delta(V_1)/\Delta(V_2) = \frac{\prod (t^{a_i}-1)}{\prod(t^{b_i}-1)}
\]
since $\wh(\ZG)\to \wh(\QG)$ is monic
 \cite[p.14]{o1}.  

Let $W$ be a finite-dimensional  $G$-representation. 
 A necessary condition for a non-linear similarity $V_1\oplus W \sim_t V_2\oplus W$ 
is the existence of a  $G$-homotopy equivalence
 $f\colon S(V_2) \to S(V_1)$ such that
$f \ast id\colon S(V_2\oplus  U)\to S(V_1\oplus  U)$ is freely 
$G$-normally cobordant to the identity map on $S(V_1\oplus  U)$,
for \emph{all} free $G${--}representations $U$
(see \pone, Section \ref{basics}).
If $V_1$ and $V_2$ satisfy this
condition, we say that $S(V_1)$ and $S(V_2)$ are \emph{$s${--}normally cobordant}. 
This  condition for non-linear similarity can
be decided by explicit congruences in the weights of $V_1$
and $V_2$ (see \cite[Thm. 1.2]{yo1}).

This quantity, $\Delta(V_1)/\Delta(V_2)$ is the basic invariant
determining non-linear similarity. It represents a unit in the group
ring $\ZG$, explicitly described for $G = \Cy{2^r}$ by Cappell and
Shaneson in
\cite[\S 1]{cs5} using a pull-back square of rings. 
To state concrete results we need to evaluate this invariant
modulo suitable indeterminacy. 

The involution $t \mapsto t^{-1}$ induces the identity on $\wh(\ZG)$,
so we get an element
\[\left\{ \Delta(V_1)/\Delta(V_2) \right\}
\in H^0(\wh(\ZG))\]
where we use $H^i(A)$ to denote the Tate cohomology
$H^i(\cy2;A)$ of $\cy2$ with coefficients in $A$.

Let $\wh(\ZG^-)$ denote the Whitehead group $\wh(\ZG)$
together with the involution induced by $t \mapsto -t^{-1}$.
Then for  $\tau(t) = \frac{\prod(t^{a_i}-1)}{\prod(t^{b_i}-1)}$,
we compute
\[\tau(t)\tau(-t)=
\frac{\prod(t^{a_i}-1)\prod((-t)^{a_i}-1)}
{\prod(t^{b_i}-1)\prod((-t)^{b_i}-1)}
= \prod\frac{(t^2)^{a_i}-1}{((t^2)^{b_i}-1)}\]
which is clearly induced from $\wh(\ZH)$.
Hence we also get a well defined
element
 \[\left\{ \Delta(V_1)/\Delta(V_2) \right\}
\in H^1(\wh(\ZG^-)/\wh(\ZH))\ .\]
This calculation takes place over the ring
 $\Lambda_{2q} =\bz [t]/(1 + t^2 + \dots + t^{4q-2})$,
 but the result holds
over
$\ZG$ via the involution{--}invariant pull-back square
\[
 \begin{matrix} \ZG  & \to &\Lambda_{2q}\cr
\downarrow &&\downarrow\cr
\bz[\cy 2]& \to & \cy{2q}[\cy 2]
\end{matrix}
\]

Consider the exact sequence of modules with involution:
\eqncount
\begin{equation}\label{okseq}
K_1(\ZH) \to K_1(\ZG)\to
 K_1(\rel{\ZH}{\ZG}) \to \widetilde K_0(\ZH)\to \widetilde K_0(\ZG)
\end{equation}
and define
$\wh(\rel{\ZH}{\ZG}) =K_1(\rel{\ZH}{\ZG})/{\scriptstyle \{\pm
G\}}$ . We then have a short exact sequence
\[0\to \wh(\ZG)/\wh(\ZH) \to \wh(\rel{\ZH}{\ZG}) \to \ck\to 0\]
where $\ck =\ker(\widetilde K_0(\ZH) \to \widetilde K_0(\ZG))$.
Such an exact sequence of $\cy 2$-modules induces a long exact
sequence in Tate cohomology.
In particular, we have a coboundary map
\[\delta\colon H^0(\ck) \to H^1(\wh(\ZG^-)/\wh(\ZH))\ .\]
Our first result deals with isotropy groups of index 2, as
is the case  for all the non{--}linear similarities
constructed in \cite{cs1}. 
\begin{thma}
Let $V_1=t^{a_1} + \dots + t^{a_k} $
and $V_2 = t^{b_1} + \dots + t^{b_k}$ be free 
$G${--}representations,
with $a_i \equiv b_i\equiv 1 \Mod{4}$.
There exists a topological similarity $V_1\oplus \rminus \sim_t
V_2\oplus \rminus$ if and only if
\begin{enumerate}
\renewcommand{\labelenumi}{(\roman{enumi})}
\item $\prod a_i \equiv \prod b_i \Mod{4q}$,
\item $\Res_H V_1 \cong \Res_H V_2$, and
\item the element
$\left\{ \Delta(V_1)/\Delta(V_2) \right\}
\in H^1(\wh(\ZG^-)/\wh(\ZH))$
is in the image of the coboundary
$\delta\colon H^0(\ck) \to H^1(\wh(\ZG^-)/\wh(\ZH))$.
\end{enumerate}
\end{thma}
\begin{remark}
The proof of this result is in Part \partone, but note that
Condition (iii) simplifies for $G$ a cyclic $2$-group
since $H^0(\ck) = 0$ in that case (see \pone, Lemma \ref{vanB}).
Theorem A  should be compared with \cite[Cor.1]{cs1},
where more explicit conditions are given for ``first-time" similarities
of this kind under the assumption
that $q$ is odd, or a $2$-power, or $4q$ is a ``tempered" number.
See also Theorem \ref{Cprime} for a more general result concerning
similarities without $\rplus$ summands.
The case $\dim V_1 = \dim V_2 = 4$ gives a reduction to
number theory for the existence of $5$-dimensional
similarities (see \pone, Remark \ref{Arem}).
\end{remark}

Our next result uses a more elaborate setting for the invariant.
Let
\[
\Phi = \begin{pmatrix} \ZH  & \to &\Zadictwo{H}\cr
\downarrow &&\downarrow\cr
\ZG & \to & \Zadictwo{G}
\end{pmatrix}
\]
and consider the exact sequence
\eqncount
\begin{equation}0 \to K_1(\rel{\ZH}{\ZG}) \to K_1(\rel{\Zadictwo{H}}{\Zadictwo{G}})
\to K_1(\Phi) \to \widetilde K_0(\rel{\ZH}{\ZG}) \to 0\ .
\end{equation}
Again we can define the Whitehead group versions by dividing
out trivial units $\{\pm G\}$, and get a double coboundary
\[\delta^2\colon H^1(\widetilde K_0(\rel{\ZH}{\ZG^{-}}))
\to H^1( \wh(\rel{\ZH}{\ZG^{-}}))\ .\]
There is a natural map
$H^1(\wh(\ZG^-)/\wh(\ZH))\to H^1( \wh(\rel{\ZH}{\ZG^{-}}))$.
We will use the same notation 
$$\{ \Delta(V_1)/\Delta(V_2) \} \in H^1( \wh(\rel{\ZH}{\ZG^{-}}))$$
 to denote the image of our Reidemeister torsion invariant.
The non-linear similarities handled by the next result have
isotropy of index $\leq 2$.
\begin{thmb}
Let $V_1=t^{a_1} + \dots + t^{a_k} $
and $V_2 = t^{b_1} + \dots + t^{b_k}$ be free $G${--}representations.
There exists a topological similarity
 $V_1\oplus \rminus\oplus \rplus \sim_t
V_2\oplus \rminus\oplus \rplus$ if and only if
\begin{enumerate}
\renewcommand{\labelenumi}{(\roman{enumi})}
\item $\prod a_i \equiv \prod b_i \Mod{4q}$,
\item $\Res_H V_1 \cong \Res_H V_2$, and
\item the element
$\left\{ \Delta(V_1)/\Delta(V_2) \right\}$
is in the image of the double coboundary
\[\delta^2\colon H^1(\widetilde K_0(\rel{\ZH}{\ZG^{-}}))
\to H^1(\wh(\rel{\ZH}{\ZG^{-}}))\ .\]
\end{enumerate}
\end{thmb}
This result can be applied to $6$-dimensional similarities.
\begin{corollary}\label{sixdim}
Let $G=\Cy{4q}$, with $q$ odd, and suppose that the
fields $\bq(\zeta_{d})$ have odd class number for
all $d \mid 4q$.
Then $G$ has no $6$-dimensional non-linear similarities.
\end{corollary}
\begin{remark} For example, the class number condition is satisfied for $q\leq 11$, but not for $q=29$.
The proof of the Corollary \ref{sixdim} is given in Section \ref{sixdiml} assuming
Theorem B, which is proved in Part \partone.
This result corrects \cite[Thm.1(i)]{cssw1}, and shows that the
computations of $\rtop{G}$ given in \cite[Thm. 2]{cssw1}
are incorrect.
\end{remark}
Our final example of the computation of  bounded transfers is suitable
for determining stable non{--}linear similarities inductively,
with only a minor assumption on the isotropy subgroups.
To state the algebraic conditions, we must again generalize the
indeterminacy for the Reidemeister torsion invariant
to include bounded $K${--}groups (see Section \ref{excision}).
In this setting $\widetilde K_0(\rel{\ZH}{\ZG}) = \widetilde K_0(\ccat{\rminus,G}{\bz})$
and $\wh(\rel{\ZH}{\ZG})=\wh(\ccat{\rminus,G}{\bz})$.
We  consider the analogous double coboundary
\[\delta^2\colon
H^1(\widetilde K_0(\ccat{W\times \rminus,G}{\bz})) \to
H^1(\wh(\ccat{W\times \rminus,G}{\bz}))\]
and note that there is a map
$\wh(\ccat{\rminus,G}{\bz})\to
\wh(\ccat{W\times \rminus,G}{\bz})$
induced by the inclusion on the control spaces. 
We will again use the same notation
$$\{\Delta(V_1)/\Delta(V_2)\}\in H^1(\wh(\ccat{W\times \rminus,G}{\bz}))$$
 for the image of the Reidemeister torsion invariant in this new domain.

\begin{thmc}
Let $V_1=t^{a_1} + \dots + t^{a_k} $
and $V_2 = t^{b_1} + \dots + t^{b_k}$ be
 free $G${--}representations.
Let $W$ be a  complex $G${--}representation
with no $\rplus$ summands.
Then there exists a topological similarity
 $V_1\oplus W \oplus \rminus\oplus \rplus \sim_t
V_2\oplus W \oplus \rminus\oplus \rplus$ if and only if
\begin{enumerate}
\renewcommand{\labelenumi}{(\roman{enumi})}
\item $S(V_1)$ is $s$-normally cobordant
 to $S(V_2)$,
\item $\Res_H (V_1\oplus W) \oplus \rplus
 \sim_t \Res_H (V_2\oplus W)\oplus \rplus$, and
\item the element
$\left\{ \Delta(V_1)/\Delta(V_2) \right\}$
is in the image of the double coboundary
\[\delta^2\colon
H^1(\widetilde K_0(\ccat{\Wmax\times \rminus,G}{\bz})) \to
H^1(\wh(\ccat{\Wmax\times \rminus,G}{\bz}))\ , \]
where $0\subseteq \Wmax \subseteq W$ is a
complex subrepresentation of real dimension $\leq 2$, with  maximal
isotropy group among the isotropy groups
of $W$ with 2-power index.
\end{enumerate}
\end{thmc}
\begin{remark}\label{unstablenormal}
 The existence of a
 similarity $V_1\oplus W\sim_t V_2\oplus W$ implies that
 $S(V_1)$ and $S(V_2)$ are  $s$-normally cobordant.
In particular, $S(V_1)$ must be freely $G$-normally cobordant to
$S(V_2)$ and this (unstable) normal invariant condition is enough
to give us a surgery problem. Crossing with $W$ defines the bounded transfer map
$$\trf{W}\colon L^h_{2k}(\ZG) \to L^h_{2k + \dim W}(\bcat{W})$$
introduced in \cite{hp1}. The vanishing of the surgery obstruction
is equivalent to the existence of a similarity (see \pone, Theorem \ref{maintest}). The computation of
the bounded transfer in  $L${-}theory leads  to condition (iii), and an
expression of the obstruction purely in terms of bounded $K$-theory. To carry
out this computation we may need  to stabilize in the free part,
and this uses the $s$-normal cobordism condition.
\end{remark}
\begin{remark} Note that $\Wmax=0$ in condition 
(iii) if $W$ has no isotropy
subgroups of $2$-power index.
Theorem C suffices to handle stable topological similarities,
but leaves out cases where $W$ has an
odd number of $\rminus$ summands (handled in  Theorem
\ref{Cprime} and the results of Section
\ref{rminustwo}). 
Simpler conditions can be given when $G = \Cy{2^r}$
 (see \pone, Section \ref{twogroups}).

The double coboundary in (iii) can also be expressed in more
``classical" terms by using the short exact sequence
\eqncount
\begin{equation} 0 \to \wh(\bCatrminus) \to 
\wh(\bcatrminus{\Wmax}) \to
K_1(\brelCatrminus{\Wmax}) \to 0       \label{ses}
\end{equation}
derived in Corollary \ref{kersimple}.
We have $K_1(\brelCatrminus{\Wmax}) = \Kminus{1}(\ZK)$, where
$K$ is the isotropy group of $\Wmax$, and
$\wh(\bCatrminus)=\wh(\rel{\ZH}{\ZG})$.
The indeterminacy in Theorem C
is then generated by the double coboundary
\[\delta^2\colon H^1(\widetilde K_0(\rel{\ZH}{\ZG^{-}}))
\to H^1(\wh(\rel{\ZH}{\ZG^{-}}))\]
used in Theorem B and the coboundary
\[\delta\colon \Coh{0}{\Kminus{1}(\ZK)} \to \Coh{1}{\wh(\rel{\ZH}{\ZG^{-}})}\]
from the Tate cohomology sequence of (\ref{ses}).

\end{remark}

Finally, we   apply these results to $\rtop{G}$. Since its rank is
known  (see  \cite{cs2} or Section \ref{ratcomp}), it remains to determine its torsion subgroup.
In Section \ref{splitting}, we will define a filtration
\eqncount
\begin{equation}R_t(G) \subseteq R_n(G) \subseteq R_h(G) \subseteq R(G)
\label{tnh}
\end{equation}
on the real representation ring
$R(G)$, inducing a filtration on
 $\rtop{G} =  R(G)/R_t(G)$.
Here the subgroup
$$R_t(G) = \{  (V_1 - V_2) \mid V_1\oplus W
\sim_t V_2\oplus W \text{ for\ some\ }W\}$$
is generated by stable topological similarity.
Note that $R(G)$ has the following nice basis:
$\{ t^i, \delta, \epsilon\,\mid\, 1\leqq i \leqq 2q-1\}$, 
where $\delta=[\rminus]$ and $\epsilon=[\rplus]$
(although  we do not have $i\equiv 1 \Mod{4}$ for all  the weights).

Let $\rfree{G} = \{ t^a \mid (a, 4q) = 1\} \subset R(G)$
 be the subgroup
generated by the free representations.
To complete the definition, we let $\rfree{\Cy 2} = \{\rminus\}$
and $\rfree{e} = \{\rplus\}$.
Then inflation and fixed sets of representations defines an isomorphism
\[R(G) = \bigoplus\limits_{K \subseteq G} \rfree{G/K}\]
 and this
direct sum splitting can be intersected with $R_t(G)$
 to define 
$\Rfree_t(G)$.
We let $\rtopfree{G} = \rfree{G}/\Rfree_t(G)$.
Since  inflation and fixed sets preserve topological similarities, we obtain an induced splitting
\[\rtopfree{G} = \bigoplus\limits_{K \subseteq G} \rtopfree{G/K}\ .\]
By induction on the order of $G$, we see that it suffices to
study the  summand $\rtopfree{G}$.

Let $\rtildefree{G} = \ker(\Res\colon\rfree{G} \to \rfree{\Godd})$,
 and then project into $\rtop{G}$ to define
\[\rtildetopfree{G} = \rtildefree{G}/\Rfree_t(G)\ .\]
In Section \ref{ratcomp} we prove that $\rtildetopfree{G}$ 
is \emph{precisely}
the torsion subgroup of $\rtopfree{G}$.
Here is a specific computation (correcting \cite[Thm. 2]{cssw1}).
\begin{thmd}
Let $G=\Cy{4q}$, with $q>1$ odd, and  suppose that the
fields $\bq(\zeta_{d})$ have odd class number for
all $d \mid 4q$.
Then $\rtildetopfree{G} = \cy{4}$ generated by
$(t-t^{1+2q})$.
\end{thmd}
For any cyclic group $G$, we use normal cobordism and homotopy equivalence to define a filtration
$$\Rfree_t(G) \subseteq	 \Rfree_n(G) \subseteq\Rfree_h(G) \subseteq\rfree{G} $$
leading by direct sum to the filtration of $R(G)$ mentioned above.
Both $\rtildefree{G}/\Rtildefree_h(G)$ and
$\Rtildefree_h(G)/\Rtildefree_n(G)$
are torsion groups  which can be explicitly
determined by congruences in the weights
 (see  Section \ref{normalinvariant} and \cite[Thm.1.2]{yo1}).
The subquotient $\Rtildefree_n(G)/\Rfree_t(G)$ always has exponent 
two (see Section \ref{sD}).

We conclude this list of sample results with a calculation of
$\rtop{G}$ for cyclic $2$-groups (see \pone\ for the proof).

\begin{thme}
Let $G=\Cy{2^r}$, with $r\geqq 4$.
Then
\[\rtildetopfree{G} = \big\langle \alpha_1, \alpha_2, \dots, 
\alpha_{r-2}, \beta_1, \beta_2, \dots,  \beta_{r-3} 
\big\rangle\]
subject to the relations $2^s \alpha_s = 0$ for $1\leqq s \leqq r-2$,
and $2^{s-1}(\alpha_s + \beta_s) =0$ for $2\leqq s \leqq r-3$, 
together with $2 (\alpha_1+\beta_1) =0$.
\end{thme}
\noindent
The generators for $r\geqq 4$ are given by the elements
\[\begin{matrix} \alpha_s = t - t^{5^{2^{r-s-2}}} & \text{and}&
\beta_s =t^5 - t^{5^{2^{r-s-2}+1}}.\end{matrix} \]
We remark that  $\rtildetopfree{\Cy{8}} = \cy 4$ generated by
$t-t^5$. 
In \pone, Theorem \ref{unstablefour}
we use this information to give a complete
topological classification of linear representations for cyclic
$2$-groups. 

\begin{acknowledgement} The authors would like to express their
appreciation to the referee for many constructive comments and suggestions. In
particular, the referee pointed out an improvement
to the statement of Theorem E. If we let $\psi_i = \alpha_{i+1}
+\beta_{i+1}$ for $1\leq i\leq r-4$, then $\rtildetopfree{\Cy{2^r}}$
is a direct sum of cyclic groups generated by $\alpha_1,\dots,\alpha_{r-2},
\beta_1$ and $\psi_1,\dots, \psi_{r-4}$
where the order of the cyclic group generated by a basis element
with subscript $i$ is $2^i$. This basis displays the group structure
explicitly.
\end{acknowledgement}

\section{ A splitting of $R_{\Top}(G)$}\label{splitting}
In this section we point out an elementary
splitting of $\rtop{G}$, and some useful filtrations.
For $G$ any finite group, we denote by $R(G)$ the \emph{real}
representation ring of $G$.
Elements in $R(G)$ can be given as formal differences
$(V_1 - V_2)$ of $G${--}representations, and
$(V_1 - V_2) \sim 0$ if and only if there exists a representation
$W$ such that $V_1\oplus W \cong V_2 \oplus W$.

Notice that for $K$ any normal subgroup of $G$,
 taking fixed sets gives a retraction of
the inflation map
\[\Inf_K\colon R(G/K) \to R(G)\]
defined by pulling back a $G/K$ representation
using the  composition with the quotient map $G \to G/K$.
More explicitly,
\[\Fix_K\colon R(G) \to R(G/K)\]
is defined by $\Fix_K (V_1 - V_2) = (V_1^K -V_2^K)$ for
each normal subgroup $K \triangleleft\, G$.
Then $\Fix_K\circ \Inf_K= id\colon R(G/K) \to R(G/K)$.
\begin{definition} A $G${--}representation $V$ is
\emph{free} if $V^K = \{ 0\}$ for all non{--}trivial normal
 subgroups $1\neq K \, \triangleleft\, G$.
\end{definition}

This is the same as the usual definition
(no non{--}identity element of $G$
fixes any non{--}zero vector) for cyclic groups.
We let
\eqncount
\begin{equation}\Rfree(G) =
\bigcap \{ \ker\Fix_K \mid 1\neq K\triangleleft\, G\}
\end{equation}
denote the subgroup of $R(G)$ such that
there is a representative
of the stable equivalence class $(V_1 - V_2)$ with
$V_1, V_2$ free representations.
\begin{proposition}There is a direct sum splitting
\[R(G) = \bigoplus\limits_{K\triangleleft\, G} \Rfree (G/K)\]
indexed by the normal subgroups $K$ in $G$.
\end{proposition}
\begin{proof}
Let $V(K)$ denote the $G${--}invariant subspace given by
the sum of all
the irreducible sub{--}representations of $V$ with kernel exactly
$K$. This is a free $G/K$ representation. 
The decomposition above is given
by mapping $(V_1-V_2)$ to the elements $(V_1(K) - V_2(K))$.
\end{proof}

Inside $R(G)$ we have the subgroup of stably topologically
similar representations
\eqncount
\begin{equation}
R_t(G) = \{  (V_1 - V_2) \mid V_1\oplus W
\sim_t V_2\oplus W \text{ for\ some\ }W\}
\end{equation}
and the quotient group is $\rtop{G}$ by definition.
We define $\Rfree_t(G) =\Rfree(G)\cap R_t(G)$.
Since $R_t(G)$ is preserved by inflation and
taking fixed sets, we obtain
\begin{corollary} There is a direct sum decomposition
\[\rtop{G} = \bigoplus\limits_{K\triangleleft\, G} \Rtopfree(G/K)\]
where the summands are the quotients
$\Rfree (G/K)\big / \Rfree_t(G/K)$.
\end{corollary}
We will also need a certain filtration of $R(G)$.
First we define
\eqncount
\begin{equation}\Rfree_h(G) = \{  (V_1 - V_2) \mid S(V_1)
\simeq_G S(V_2) \text{ \ for\ $V_1$\ and $V_2$\ free}\}
\end{equation}
where $\simeq_G$ denotes  $G${--}homotopy equivalence.
This is a subgroup of $\Rfree(G)$, in fact a sub{--}Mackey functor
since it has induction and restriction for subgroups of $G$.
We define 
\eqncount
\begin{equation}
R_h(G) = \bigoplus\limits_{K\triangleleft\, G} \Rfree_h(G/K)
\end{equation}
If there exists a $G${--}homotopy equivalence
 $f\colon S(V_1) \to S(V_2) $ such that
\[S(V_1 \oplus U)=S(V_1){\ast} S( U)
\RA{\ f{\ast} 1\ } S(V_2){\ast} S( U) = S(V_2\oplus U)
\]
is freely $G${--}normally cobordant to the identity for all free
$G${--}representations $U$, then we say that
$S(V_1)$ and $S(V_2)$ are $s$-normally cobordant, and we
write $S(V_1) \bumpeq_G S(V_2)$. Define
\eqncount
\begin{equation}
\Rfree_n(G) = \{ \alpha\in \Rfree_h(G)  \mid \exists V_1, V_2 \text{\ with\ }\alpha=(V_1-V_2)
 \text{ \ and\ } S(V_1)
\bumpeq_G S(V_2)\}
\end{equation}
and note that $\Rfree_n(G)$ is also a subgroup of $\Rfree(G)$.
Indeed, if $(V_1 - V_2)$ and $(V'_1 - V'_2)$
are in $\Rfree_n(G)$, there exist $G$-homotopy equivalences
 $f\colon S(V_1) \to S(V_2) $
and $f'\colon S(V'_1) \to S(V'_2) $
with $f{\ast} 1$ and $f'{\ast} 1$
normally cobordant to the identity
under any stabilization.
But $f{\ast} f' \simeq_G (1{\ast} f')\circ (f{\ast} 1)$, so we
just glue together the normal cobordisms
for $f{\ast} 1$ (after stabilizing by $U=V'_1$) and for
$1{\ast} f'$ (after stabilizing by $U=V_2$) along
the common boundary $id\colon S(V_2 \oplus V'_1)$.
 As above, we define
 \eqncount
 \begin{equation}
R_n(G) = \bigoplus\limits_{K\triangleleft\, G} \Rfree_n(G/K)
\end{equation}
  Since $\Rfree_t(G) \subseteq \Rfree_n(G)$, we have defined a filtration
\eqncount
\begin{equation}
R_t(G) \subseteq R_n(G) \subseteq R_h(G) \subseteq R(G)\label{filt}
\end{equation}
of $R(G)$, natural with respect to  
restriction of representations. All the terms except possibly
$R_n(G)$ are also natural with respect to induction
of representations.

\begin{remark} It follows from the proof of \cite[3.1]{yo1} that 
$S(V_1\oplus U_0)$ is $s$-normally cobordant to
$S(V_2\oplus U_0)$, for some free $G$-representation $U_0$, 
if and only if $S(V_1)$ is $s$-normally cobordant to
$S(V_2)$. It follows that we could have used the latter condition
to define $\Rfree_n(G)$.
\end{remark}

\section{ A rational computation}
\label{ratcomp}
In this section we use \pone,  Theorem \ref{maintest} and 
the splitting of the last
section to describe the torsion subgroup of $\rtopfree{G}$.
We also give a new proof of Cappell and Shaneson's result
computing $\rtop{G}\otimes \bq$  for all finite groups $G$.

First we consider cyclic groups.
 Let $G=\Cy{2^rq}$ be a cyclic group, where $q\geqq 1$ is odd.
The Odd Order Theorem \cite{hsp1},
\cite{mr1} gives $R(G) = \rtop{G}$ if $r \leqq 1$, and we
recall the definition
\[\Rtildetopfree(G) = \Rtildefree(G)/\Rfree_t(G)\]
where
$\Rtildefree(G) = \ker\{ \Res\colon \rfree{G} \to
\rfree{\Godd}\}$.
Here $\rtildefree{\Cy{2q}} = \rtildefree{\Cy{q}} = 0$,
for $q\geqq 1$ odd.

\begin{theorem}\label{threeone}
For $G$ cyclic,
the kernel and cokernel of
\[\Res\colon \rtopfree{G} \to \rtopfree{\Godd}=\rfree{\Godd}\]
are $2${--}primary torsion groups.
\end{theorem}
\begin{corollary}\label{cyclicranks}
Let $G=\Cy{2^rq}$, $q$ odd, be a finite cyclic group.
\begin{enumerate}
\renewcommand{\labelenumi}{(\roman{enumi})}
\item The torsion subgroup of $\rtopfree{G}$
is $\rtildetopfree{G}$.
\item The rank of $\rtopfree{G}$ is $\varphi(q)/2$
(resp. $1$) if $q>1$ (resp.  $q=1$).
\item We have the formula
\[\text{rank} (\rtop{G}) =
\begin{cases} (r+1)\bigl\{
 \sum_{1 \not=d \mid q} \varphi(d)/2 + 1\bigr\},& \text{for $q>1$}\cr
(r+1)& \text{for $q=1$}.
\end{cases}
\]
where $\varphi(d)$ is the Euler function.
\end{enumerate}
\end{corollary}
\begin{proof}
The first part follows directly from
Theorem \ref{threeone}. For parts (ii) and (iii)
 count the free representations of $\Godd$.
\end{proof}
Now let $\Kodd$ be the composite of all subfields of $\br$
of the form $\bq (\zeta+\zeta^{-1})$ where
$\zeta \in \bc$ is an odd root of unity.
\begin{lemma}\label{threetwo}
For $G$ cyclic,
the composition of the natural map $R_\Kodd(G) \to R(G)$ and
$\Res\colon R(G) \to R(\Godd)$ induces
a $p${--}local isomorphism $\Rfree_\Kodd(G) \to \Rfree(\Godd)$,
for any odd prime $p$.
\end{lemma}
\begin{proof}
According to a result of Brauer \cite[Thm.24]{se1}
 $R_\Kodd(\Godd) = R(\Godd)$, and any representation
of $G$ can be realized over the field $\bq (\zeta_{|G|})$.
In addition,  the restriction map
 \[\Res\colon R_\Kodd(G) \to R_\Kodd(\Godd)\]
is a $p${--}local surjection
(since $\Res_{\Godd}\circ \Ind_{\Godd}$ is just multiplication by $2^r$).
But the rank of $R_\Kodd(G)$  given in \cite[12.4]{se1} equals the rank of
$R(\Godd)$, so we are done.
\end{proof}
\begin{corollary}\label{threethree}
For $G$ cyclic and any odd prime $p$,
the natural map  $R_\Kodd(G) \to \rtop{G}$ induces
a $p${--}local isomorphism.
\end{corollary}
\begin{proof}
This follows from the Lemma \ref{threetwo} and
Theorem \ref{threeone}.
\end{proof}
In \cite{cs2}, Cappell and Shaneson obtained
the following result by a different argument. It computes the rank of
$\rtop{G}\otimes \bq$ for any finite group $G$.
\begin{theorem}
Let $\Kodd$ be the composite of all subfields of $\bc$
of the form $\bq (\zeta+\zeta^{-1})$ where
$\zeta \in \bc$ is an odd root of unity.
Then for any finite group $G$,  the natural map $R_\Kodd(G) \to R(G)$ induces
an isomorphism $R_\Kodd(G)\otimes \bz _{(p)}
\cong \rtop{G}\otimes \bz _{(p)}$ for any
odd prime $p$.
\end{theorem}

\begin{proof}
The result holds for cyclic groups $G$ by Corollary \ref{threethree},
and we apply induction theory to handle general finite groups.

First we observe that the Mackey functors
$R(G)\otimes \bz _{(p)} $ and $R_\Kodd(G)\otimes \bz _{(p)}$
are  generated by induction from
$p${--}elementary subgroups (see
\cite[Thm.27]{se1}  and note that
$\Gamma_\Kodd = \{\pm 1\}$).
Therefore,  by \cite[11.2]{o1}
they are $p${--}elementary computable in the sense
of Dress induction theory \cite{dre1}.
We may therefore assume that $G$ is $p${--}elementary.

Now if $G$ is $p${--}elementary, it is a product of a
$p${--}group and a cyclic group prime to $p$.
Any irreducible complex representation of $G$ is then
induced from a linear representation on a subgroup
 \cite[Thm.16]{se1}.

It follows that $\Rtop(G)$ is generated by  generalized induction
(i.e. inflation followed by induction) from cyclic subquotients.
Consider the following commutative diagram
\[
\xymatrix{R_\Kodd(G) \ar[r]\ar@<1ex>@{  >->}[d]^{\Res}&
R(G)\ar[r]\ar@<1ex>[d]&
\rtop{G}\ar@<1ex>[d]\\
{\bigoplus\limits_{C}}
 R_\Kodd(C) \ar@<1ex>[u]^{\Ind}\ar@<1ex>[r]&
{\bigoplus\limits_{C}}
R(C)\ar@<1ex>[u] \ar@<1ex>[r]&
{\bigoplus\limits_{C}}
 \rtop{C}\ar@<1ex>@{>>}[u]
}
\]
It follows from
 Corollary \ref{threethree}  that the top composite
 map $R_\Kodd(G)\otimes \bz _{(p)} \to \rtop{G})\otimes \bz _{(p)}$
is surjective.

The sum of the (ordinary) restriction maps
to cyclic subgroups induces a rational injection on $R_\Kodd(G)$
(see \cite[2.5,2.10]{swe1}).
Since $R_\Kodd(G)$ is torsion{--}free, it follows again from
Corollary \ref{threethree} that the map
 $R_\Kodd(G)\otimes \bz _{(p)} \to \rtop{G})\otimes \bz _{(p)}$
is injective.
\end{proof}

\begin{proof}[The proof of Theorem \ref{threeone}]
Since an $\Kodd${--}representation is free if and only if it
is the sum of Galois conjugates of free $G${--}representations,
we can decompose $R_\Kodd(G)$ as in Section \ref{splitting}, and
conclude that $\Rfree_\Kodd(\Godd) = \Rfree(\Godd)$.
 It remains to show that $\Rtildetopfree(G) $ is
a torsion group with 2{--}primary exponent.
For this we use the filtration of \S \ref{splitting}.

For $\Rtildefree(G)/\Rtildefree_h(G) $ this is easy
since the $k${--}invariant gives a homomorphism (via joins
of free  $G${--}spheres) to $\Units{(\cy{2^r})}$ and this is
a 2{--}group.
The next quotient is also 2{--}primary torsion, by results
of \cite{yo1}: a sufficiently large 2{--}power join of a $G${--}homotopy
equivalence between two  free $G${--}spheres,
which are linearly equivalent over $\Godd$, becomes $s${--}normally
cobordant to the identity.
The point is that the normal invariant is detected by a
finite number of 2{--}power
congruences
conditions among the Hirzebruch $L${--}classes of the tangent bundles
of the lens spaces, and this can be satisfied after sufficiently
many joins.

Finally, the last quotient $\Rtildefree_n(G)/ \Rtildefree_t(G) $
is shown to be $2${--}primary torsion in the next proposition.
\end{proof}
\begin{proposition}
Let $G=\Cy{2^rq}$, $q$ odd,  and assume
\[\sigma \in \ker(\Res \colon  L^h_{2k}(\ZG) \to L^h_{2k}(\bz \Godd)).\]
Then there exists a complex representation
$W$ with $W^G =0$ such  that $\trf{W}(2^r\sigma)=0$.
\end{proposition}
\begin{proof}
We will take $W = \rminus \oplus \rminus \oplus W_0$,
where $W_0$ is the sum of all the irreducible $2${--}dimensional
representations of $G$ with isotropy of $2${--}power index.
Note that the $\rminus$-transfer is just the compact
$I_{-}$ transfer of one-sided codimension 1 surgery followed by
adding rays to infinity, so whenever the $I_{-}$transfer is $0$, the
$\rminus$-transfer will have to be $0$.
This was discussed in more detail
in \pone, Section \ref{rmin}.

\smallskip

\noindent\textbf {Step 1:}
If $G$ has odd order, there is nothing to prove.
Otherwise, let $H\subset G$ be of index 2.
If $\Res_H(\sigma)=0$, then
$\trfrminus(\sigma) \in
\Image(L_1^h(\ZG^-) \to L_1^h(\rel{\ZH}{\ZG^-}))$.
But  $L_1^h(\ZG^-)$ has exponent 2
\cite[12.3]{ht1}, so $\trfrminus(2\sigma) =0$.
Then take $W=\rminus\oplus\rminus$.
Note that this case applies to $G=\Cy{2q}$,
so we can always get started.

\smallskip

\noindent\textbf{Step 2:}
We may assume that $r\geqq 2$.
If $\Res_H(\sigma)\ne 0$ note that
\[\Res_H(2\sigma -\Ind_H\Res_H(\sigma))
= 2\Res_H(\sigma) - 2  \Res_H(\sigma) = 0.\]
By induction  $\Res_H W_0$
works for $\Res_H(\sigma)$: say
\[\trf{\Res_H W_0}(2^{r-1}\Res_H(\sigma))=0\]
 and $W_0^H=0$.
 Let $\dim W_0=m$ and consider the commutative diagram
\[
\xymatrix{ L_{2k}(\ZG) \ar[r]^
{\Res} \ar[d]_{\trf{W_0}}&L_{2k}(\ZH)\ar[r]^{\Ind}
\ar[d]_{\trf{\Res_H W_0}}&L_{2k}(\ZG)\ar[d]_{\trf{W_0}}\\
L_{2k+m}(\mathcal C_{W_0,G}(\bz))\ar[r]^{\Res}& L_{2k+m}(\mathcal
C_{\Res W_0,H}(\bz))\ar[r]^{\Ind}&L_{2k+m}(\mathcal C_{W_0,G}(\bz)).
}\]
 From this we get
\[2^{r-1}\cdot\trf{W_0}(\Ind_H\Res_H(\sigma)) = 0\ .\]
The first step implies that $\trf{\Rminus{2}}(2\sigma_1) = 0$,
where $\sigma_1 = 2\sigma - \Ind_H\Res_H(\sigma)$.
Let $W= \rminus\oplus\rminus\oplus W_0$, so that
we have $W$ complex and $W^G=0$.
Note that $\trf{W}= \trf{\Rminus{2}}\circ\trf{W_0}=
\trf{W_0}\circ\trf{\Rminus{2}}$.
But
\[2^r\cdot \trf{W}(\sigma)
= 2^{r-1}\trf{W}(2\sigma - \Ind_H\Res_H(\sigma)) +
2^{r-1}\trf{W}(\Ind_H\Res_H(\sigma))\]
and both terms vanish
(because $r\geqq 2$ and by the property of $W_0$
respectively).
\end{proof}
A similar argument to that in Step 2 above gives:
\begin{proposition}\label{trftorsion}
If $\Res_H(\trf{W}(x))=0$ for $x \in L_0^h(\ZG)$,
then $4\cdot\trf{W\times\rminus}(x) =0$.
\end{proposition}
\begin{proof}
Since
\[4\cdot \trf{W\times\rminus}(\sigma)
= 2\cdot\trf{W}\trfrminus(2\sigma - \Ind_H\Res_H(\sigma)) +
2\cdot\trfrminus\trf{W}(\Ind_H\Res_H(\sigma))\]
we conclude as above that both terms vanish.
\end{proof}

\section{ Excision in bounded surgery theory}
\label{excision}

A small additive category with involution $\ca$ is a small additive
category together with a contravariant endofunctor $*$ such that
$*^2= 1_\ca$.
 Ranicki defines algebraic $L${--}theory $L^h_*(\ca)$ for such
categories and corresponding spectra
$\BL^h(\ca)$ with $L_*^h(\ca)=\pi_*(\BL^h(\ca))$ \cite{ra5}.
The obstruction
groups for bounded surgery are obtained this way for appropriately
chosen additive categories.
We shall also need a simple version of such
groups.
 For this, the additive category must come equipped with a
system of stable isomorphisms and a subgroup $s\subset K_1(\ca)$,
such that any composition resulting in an automorphism defines an
element in $s$.
 The point here is that whenever two objects are stably
isomorphic, there is a canonically chosen stable isomorphism,
canonical up to automorphisms defining elements of $s$.
 In this situation
Ranicki refines the definition of $\BL^h(\ca)$ to give the simple
 $L${--}theory
spectrum $\BL^s(\ca)$, by requiring appropriate isomorphisms to give
elements of $K_1(\ca)$ belonging to the subgroup $s$.
More generally, we also get the
$\BL^ {a}(\ca)$-spectra for any involution invariant subgroup
$a$ with $s \subset a \subset K_1(\ca)$, coinciding
with $\BL^ h(\ca)$ when $a = K_1(\ca)$.
\begin{example} Let $\ca$ be the category of free $\ZG$-modules
with a
$G$-invariant $\bz$-basis, and $\ZG$-module morphisms.
Two objects are stably isomorphic if and only if they have the same rank.
The preferred isomorphisms are chosen to be the ones sending a
$\bz$-basis to a
$\bz$-basis, so automorphisms define elements of $\{\pm G\}\subset
K_1(\ZG)$.
In this situation one obtains Wall's  $L^s$-groups.
\end{example}
The theory of projective $L${--}groups fits into the scheme
as follows: one defines  $\BL^p(\ca) = \BL^h(\ca^{\wedge})$, where
$\ca^\wedge$ is the idempotent completion of $\ca$.
The objects of $\ca^\wedge$ are pairs $(A,p)$ with $A$ an object of
$\ca$ and $p^2=p$.
The morphisms $\phi\colon(A,p) \to (B,q)$ are the  $\ca$-morphisms
$\phi\colon A\to B$ with $q\phi p = \phi$.
Again it is possible to ``partially'' complete $\ca$.
If $K_0(\ca)\subset \ck\subset K_0(\ca^\wedge)$ is an
involution invariant subgroup , we define $\ca^{\wedge k}$ to be the full
subcategory of $\ca^\wedge$ with objects defining elements of
$\ck\subset K_0(\ca^\wedge)$.
This way we may define $\BL^k(\ca) =
\BL(\ca^{\wedge k})$. Similarly to the above, for $\ck=K_0(\ca)\subset
K_0(\ca^\wedge)$, $L_*^k(\ca)$ is naturally isomorphic to $L^h_*(\ca)$.
The quotient $\widetilde K_0(\ca) = K_0(\ca^\wedge)/K_0(\ca)$
is called the reduced projective class group of $\ca$.
\begin{example} If $\ca$ is the category of free $\ZG$-modules
then
$\ca^\wedge$ is isomorphic to the category of projective
$\ZG$-modules and the $L^p_*(\ca)$ are Novikov's original $L^p$-groups.
\end{example}
Suppose $M$ is a metric space with the finite group $G$
acting by isometries, $R$ a ring with involution. In \cite[3.4]{hp1}
we defined an additive
category ${\mathcal G}_{M,G}\,(R)$
 with involution as follows:
\begin{definition} An object A is a  left $R(G)$-module together
with a map $f\colon A \longrightarrow F(M)$, where $F(M)$ is the set of finite
subsets of $M$, satisfying
\begin{enumerate}
\renewcommand{\labelenumi}{(\roman{enumi})}
\item $f$ is
$G$-equivariant.
\item $A_x = \{ a \in A \mid f(a) \subseteq \{ x \} \}$ is a finitely generated
free sub $R$-module
\item As an $R$-module $A = \bigoplus_{x \in M} A_x$
\item $f(a+b) \subseteq f(a) \cup f(b)$
\item For each ball $B \subset M$, $\{x \in B\mid A_x \ne 0\}$ is finite.
\end{enumerate} A morphism $\phi\colon A \longrightarrow B$ is a morphism of
$RG$-modules, satisfying the following condition: there exists $k$ so
that the components $\phi^m_n \colon A_m \longrightarrow B_n$ (which are
$R$-module morphisms) are zero when $d(m,n) > k$. 
 Then ${\mathcal G}_{M,G}\,(R)$
 is an additive category in an obvious way.
 The full subcategory of ${\mathcal G}_{M,G}\,(R)$, for which all the
 object modules are required to be \emph{free} left $RG$-modules, is denoted
 $\ccat{M,G}{R}$.

Given an object $A$, an $R$-module homomorphism $\phi \colon A
\longrightarrow R$ is said to be locally finite if the set of $x \in M$ for
which $\phi(A_x) \neq 0$ is finite. We define $A^* = \text{Hom}_R^{lf}(A,R)$,
the set of locally finite $R$-homomorphisms.
We want to make $*$ a functor from ${\mathcal G}_{M,G}\,(R)$
  to itself to make ${\mathcal G}_{M,G}\,(R)$ a category
with involution.
 We define $f^* \colon A^* \longrightarrow FM$ by $f^*(\phi) =
\{x\mid\phi(A_x) \neq 0\}$ which is finite by assumption. $A^*$ has an
obvious right action of $G$ turning it into a right $RG$ module given by
$\phi g(a) = \phi(ga)$, and $f^*$ is equivariant with respect to the right action on
$M$ given by $xg = g^{-1}x$.
To make $*$ an endofunctor of ${\mathcal G}_{M,G}\,(R)$
we need to replace the right action by a left action.
We do this by the standard way in surgery theory by  letting $g$ act on the
left by letting $g^{-1}$ act on the right.
In the unoriented case, given a
homomorphism $w\colon   G \longrightarrow
\{\pm 1\}$, we let $g$ act on the left of $A^*$ by $w(g)\cdot g^{-1}$ on
the right. The involution $*$ induces a functor on the subcategory  
$\ccat{M,G}{R}$, so that $\ccat{M,G}{R}$
 is also a category with involution.
\end{definition}
\begin{example}\label{standard}Let $\rho_W\colon G \to O(W)$ be an orthogonal
action of $G$ on a finite dimensional real vector space $W$.
We take $M=W$ with the action through $\rho_W$, and orientation
character $\det(\rho_W)$.
This will be called the \emph{standard orientation} on $\bcat{W}$.
\end{example}

\begin{remark} We will need to find a system of stable
isomorphisms for the category $\ccat{M,G}{R}$ to be able to do
simple $L${--}theory.
To do this we choose a point $x$ in each $G$-orbit,
and an $RG_x$-basis for $A_x$, where $G_x$ is the isotropy subgroup
of $x$.
We then extend that by equivariance to an $R$-basis of the
module.
Having a basis allows defining an isomorphism in the usual
fashion.
In each case we need to describe the indeterminacy in the
choices coming from the choice of $R$-basis and points in the orbit.
For our particular choices of $M$ it will be easy to determine the subgroup $s$,
so we will not formulate a general statement.
\end{remark}

We will study the $L${--}theory of the categories $\ccat{M,G}{R}$ using
excision.
Let $N$ be a $G$-invariant metric subspace of $M$.
 Denoting
$\ccat{M,G}{R}$ by $\cu$, let $\ca$ be the full subcategory on modules
$A$ so that $A_x=0$ except for $x$ in some bounded neighborhood of
$N$.
The category $\ca$ is isomorphic to $\ccat{N,G}{R}$ in an obvious
fashion.
The quotient category $\cu/\ca$, which we shall denote by
$\Ccat{M,G}{N}{R}$, has the same objects as $\cu$ but two morphisms
are identified if the difference factors through $\ca$, or
 in other words, if they  differ in a bounded neighborhood of $N$.
This is
a typical example of an additive category $\cu$ which is $\ca$-filtered
in the sense of Karoubi.
 We recall the definition.
\begin{definition} Let $\ca$ be a full subcategory of an additive
category $\cu$.
Denote objects of $\ca$ by the letters A through F and
objects of $\cu$ by the letters U through W.
We say that $\cu$ is
$\ca$-filtered, if every object $U$ has a family of decompositions
$\{U=E_\alpha\oplus U_\alpha\}$, so that
\begin{enumerate}
\renewcommand{\labelenumi}{(\roman{enumi})}
\item For each $U$, the decompositions form a filtered poset under the
partial order that $E_\alpha\oplus U_\alpha \le E_\beta\oplus U_\beta$,
whenever $U_\beta \subseteq U_\alpha$ and $E_\alpha \subseteq
E_\beta$.
\item Every map $A \to U$, factors $A \to E_\alpha \to E_\alpha\oplus
U_\alpha = U$ for some $\alpha$.
\item Every map $U\to A$ factors $U=E_\alpha \oplus U_\alpha
\to E_\alpha
\to A$ for some $\alpha$.
\item For each $U,V$ the filtration on $U\oplus V$ is equivalent to the
sum of filtrations $\{U=E_\alpha \oplus U_\alpha\}$ and $\{ V = F_\beta
\oplus V_\beta\}$ i.e. to $U\oplus V = (E_\alpha\oplus F_\beta) \oplus
(U_\alpha\oplus V_\beta)$
\end{enumerate}
\end{definition}
The main excision results were proved in \cite{pw1},
\cite{cap1}, \cite{cp1}, \cite{ra7}.
 We give a slight generalization of the $L${--}theory results.
Let $\BK $ denote  the Quillen $K${--}theory spectrum, and
$\BK ^{-\infty}$   its non-connective delooping
(with the $\Kminus{i}$-groups as homotopy groups).
\begin{theorem}\label{ua}
Let $\cu$ be an $\ca$-filtered additive category with
involution.
Consider the map
$i\colon  K_0(\ca^\wedge) \to K_0(\cu^\wedge)$
induced by inclusion, and
 let $k= i^{-1}(K_0(\cu))$.
There are fibrations of spectra
\[\BK (\ca^{\wedge k})\to \BK (\cu) \to \BK (\cu/\ca)\]
and
\[\BK ^{-\infty}(\ca)\to \BK ^{-\infty}(\cu)\to \BK ^{-\infty}(\cu/\ca)\]
If $\cu$ and $\ca$ admit compatible involutions there is a
fibration of spectra
 \[\BL^k(\ca)\to \BL^h(\cu)\to \BL^h(\cu/\ca)\ .\]
More generally, if
\begin{enumerate}
\renewcommand{\labelenumi}{(\roman{enumi})}
\item $a\subset K_i(\ca)$,$b\subset K_i(\cu^\wedge)$, and
$c\subset K_i((\cu/\ca)^\wedge)$, for $i\leqq 1$,
\item  $a= i^{-1}(b)$ and $b\to c$ is onto,
\item  $a$, $b$, and $c$ contain $K_0(\ca)$, $K_0(\cu)$ and
$K_0(\cu/\ca)$ respectively, if $i=0$, and
\item  $a$, $b$ and $c$
contain the indeterminacy subgroup given by the system of stable
isomorphisms in the case $i=1$,
\end{enumerate}
then we have a fibration of spectra
\[\BL^a(\ca)\to \BL^b(\cu) \to \BL^c(\cu/\ca)\]
\end{theorem}
\begin{proof} The $K${--}theory statements are implicitly contained in
 \cite{pw1}.
 A simpler, more modern proof and explicit statements are given in
\cite{cap1}.
The first $L${--}theory statement was proved in \cite{cp1},
and the other $L${--}theory statements follow by the following argument:
we have an exact sequence
\[K_0(\ca^{\wedge k})\to K_0(\cu)\to K_0(\cu/\ca) \to 0,\]
where  the map from $K_0(\cu)\to K_0(\cu/\ca)$ is
onto because the categories have the same objects.
Letting $I$ denote
the image of $K_1(\cu/\ca)$ in $K_0(\ca^{\wedge k})$,
we consider the diagram of short exact sequences:
\[
\xymatrix{ 0\ar[r]& I \ar@{=}[d]\ar[r]& K_0(\ca^{\wedge
k})\ar[d]\ar[r]& K_0(\cu)\ar[r]\ar[d]& K_0(\cu/\ca)\ar[r]\ar@{=}[d]& 0\\
0\ar[r]& I\ar[r]\ar@{=}[d]& a\ar[r]\ar@{=}[d]& b'\ar[d]\ar[r]&
K_0(\cu/\ca)\ar[d]\ar[r]& 0\\ 0\ar[r]& I\ar[r]\ar@{=}[d]& a\ar[d]\ar[r]&
b\ar[r]\ar[d]& c\ar[r]\ar[d]& 0\\ 0\ar[r]& I\ar[r]&
K_0(\ca^\wedge)\ar[r]& K_0(\cu^\wedge)\ar[r]& K_0((\cu/\ca)^\wedge) }
\]
 The vertical arrows are either equalities or inclusions.
We define $b'$ simultaneously as the pullback of $a/I\to b\to c$
and as the pushout of
$0\to K_0(\ca^{\wedge k})/I\to K_0(\cu)\to K_0(\cu/\ca)\to 0$.
We have
\[a/K_0(\ca^{\wedge k}) \cong  b'/K_0(\cu),\]
so using the Ranicki{--}Rothenberg
 fibrations of spectra \cite{ra7}
\[
\xymatrix{
\BL^k(\ca)\ar[r]\ar[d]&
 \BL^a(\ca)\ar[r]\ar[d]&
  {\mathbb H}(a/K_0(\ca^{\wedge k}))\ar[d]\\
\BL^h(\cu)\ar[r]& \BL^{b'}(\cu) \ar[r]& {\mathbb H} (b'/K_0(\cu))
}
\]
we get a fibration
\[\BL^a(\ca)\to \BL^{b'}(\cu)\to \BL^h(\cu/\ca).\]
We now repeat this
argument using the isomorphisms $b/b'\cong c/K_0(\cu/\ca)$
 to obtain the desired fibration of spectra.
 Since $L^h$-groups may be understood as simple
$L$-groups with all of $K_1$ as allowed torsions,
the above  bootstrapping argument
extends to fibrations of the $L$-spectra stated, using the isomorphism
\[K_1(\cu/\ca)/\ker(\bd) \cong \text{image}(\bd)\]
 where
$\bd$  is the boundary map $\bd \colon   K_1(\cu/\ca) \to
K_0(\ca^\wedge)$.
\end{proof}

In Section \ref{rminustwo} we need to use bounded surgery groups with
geometric anti{--}structure generalizing
the definition of \pone, Section \ref{rmin} (see \cite{hp1}).
The new ingredient is
 a counterpart  to the automorphism $\theta\colon H \to H$
at the metric space level.

Let  $\theta_H \colon H\longrightarrow H$ be a group automorphism
so that the data $(H, \theta_H, w, b)$ gives a geometric
anti{--}structure on $RH$.
 Let  $\theta_M \colon M \longrightarrow M$ be an isometry
 with the properties
 $\theta_M(g \cdot m) = \theta_H(g) \cdot \theta_M(m)$,
$\theta_M^2(m) = bm$, and  $\theta^2_H(g) = bgb^{-1}$.

Given an object $A \in {\mathcal G}_{M,H}\,(R)$,  we have
the functor  $*$  from
 ${\mathcal G}_{M,H}\,(R)$ to itself so that ${\mathcal G}_{M,H}\,(R)$
 is a category with involution.
We may then twist the involution $*$ by composing with the
functor sending $(A,f)$ to
$(A^{\theta},f^{\theta})$ where $A^{\theta}$ is the
 same $R$-module, but $g$ acts on the left by multiplication by
$\theta(g)$ and $f^{\theta} =\theta^{-1}_M \cdot f$.
 This defines the bounded anti{--}structure on
 ${\mathcal G}_{M,H}\,(R)$ and on the subcategory
$\ccat{M,H}{R}$ of free $RH$ modules.

\begin{example}\label{bdLNseq}
Bounded geometric anti{--}structures arise geometrically
as above.
The bounded $\rminus$ transfer sits in the long
 exact sequence
 \begin{multline}
 LN_n(\bcatrminus{W},w\phi) \to
L^h_n(\bcat{W},w) \to
L^h_{n+1}(\bcatrminus{W},w\phi)\\
 \to LN_{n-1}(\bcatrminus{W},w\phi)\to
 L^h_{n-1}(\bcat{W},w) \to\dots
\end{multline}
where $w=\det(\rho_W)$ is the standard orientation
(see Example \ref{standard}).
The bounded $LN${--}group
\[LN_n(\bcatrminus{W},w\phi)
\cong L_{n}(\ccat{W,H}{\bz},\alpha,u)\]
where $\theta_W(x) = t\cdot x$ and $\theta_H(h) = tht^{-1}$
for a fixed $t \in G-H$.

Conversely, given a bounded geometric antistructure
$(\theta_H,\theta_M,b,w)$,
we can define $G=\langle H,t\mid t^{-1}ht=\theta_H(h),t^2=b\rangle$
and $t\cdot m = \theta_M(m)$.
Then $\ccat{M,G}{R}$ induces $(\theta_H,\theta_M,b,w)$
as above, showing that all geometric antistructures
arise by twisting and restricting to an index two subgroup.
\end{example}
The $L$-theory of these bounded geometric anti{--}structures
also has a useful vanishing property which we now wish to
formulate.
We first give a basic construction.
\begin{definition}
If $\ca$ is an additive category,
then the \emph{opposite category}  $\ca^{op}$
is the category with the same objects as $\ca$
 but $\hom_{\ca^{op}}(A,B)= \hom_{\ca}(B,A)$.
The product category $\ca\times\ca^{op}$
is an additive category with involution given by
 $\ast\colon (A,B)=(B,A)$ on objects and
$\ast\colon (\alpha,\beta)=(\beta,\alpha)$
on morphisms.
\end{definition}

Clearly $K_i(\ca^{op})=K_i(\ca)$, so we can identify
$K_i(\ca \times \ca^{op}) = K_i(\ca) \times K_i(\ca)$.

\begin{lemma}\label{Rop}
 Let $b \subseteq K_i(\ca)$  for some $i\leqq 1$, and
$q=b\times b\subseteq K_i(\ca\times\ca^{op})$.
 Then  $L_{n}^q(\ca\times\ca^{op})=0$ for all $n$.
\end{lemma}
\begin{proof}
Let ${\mathcal P}(\ca)$ denote the category with the
same objects as $\ca$, but
with morphisms given by $\ca${--}isomorphisms.
Then it suffices to prove that the quadratic category
${\mathcal Q}(\ca \times \ca^{op}) \simeq {\mathcal P}(\ca)$
via the hyperbolic map (see \cite[p.122]{w8}).
This shows that $L^h_*(\ca \times \ca^{op}) =0$ and
other decorations follow trivially from the
Ranicki Rothenberg exact sequences
(note that the Tate cohomology $H^*(q)=0$).
The result for
 lower $L${--}groups follows by replacing
 $\ca$ by $\cc_{\br}(\ca)$.

Suppose $(\nu_1,\nu_2)\colon (A,B)\to (B,A)$ is a non{--}singular
quadratic form representing an element in
 ${\mathcal Q}(\ca \times \ca^{op})$.
 This means that the bilinearization $\nu_1 +\nu_2$
is an isomorphism, and we are allowed to change
$(\nu_1,\nu_2)$ by terms of the
form $(\alpha,\beta)-(\beta,\alpha)$.
We have
\[
(\nu_1,\nu_2) + (\nu_2,0) -(0,\nu_2) = (\nu_1 +\nu_2,0)
\]
and the right hand side is a hyperbolic form.
\end{proof}
We encounter the $\ca\times \ca^{op}$ situation in the
following setting:
\begin{example}
Let $M = M_1\cup M_2$ be a metric space
given as the union  of two sub{--}metric spaces $M_1$
and $M_2$, where we denote $M_1\cap M_2$ by $N$.
Suppose that
\begin{enumerate}
\renewcommand{\labelenumi}{(\roman{enumi})}
\item  $G$ acts by isometries on $M$, such that
each $g \in G$ preserves or switches $M_1$ and
$M_2$ in this decomposition,
\item  $H=\{ g\in G\mid g(M_1)=M_1\}$
is an index two subgroup of $G$,
\item for every $k>0$
 there exists an $l>0$ such that, if $x \in M_1$ (resp. $x \in M_2$)
with $ d(x,N) > l$, then
$d(x,M_2) > k$ (resp. $d(x,M_1) > k$).
\end{enumerate}
The category $\cc_{M,H}(R)$ has a
bounded geometric antistructure $(\alpha,u)$ given by
$(\theta_H,\theta_M,b,w)$ as in Example \ref{bdLNseq}, with
$\theta_M(m) = t\cdot m$ and $\theta_H(h)=tht^{-1}$
 for a fixed $t \in G -H$.
Next, observe that the category
 \[\cc_{M,H}^{>N}(R) = \cc_{M_1,H}^{>N}(R) \times \cc_{M_2,H}^{>N}(R)\]
because of our separation condition (iii).
Moreover, the functor
 $T\colon \cc_{M_1,H}^{>N}(R) \to \cc_{M_2,H}^{>N}(R)^{op}$
defined by $T(A,f) = (A^*, \theta^{-1}_M\circ f^*)$ on objects
and $T(\phi) = \phi^*$ on morphisms is an equivalence of categories.
We are thus in the $\ca\times \ca^{op}$ situation described above and
$L^h_{\ast}(\cc_{M,H}^{>N}(R),\alpha,u)=0$ by Lemma \ref{Rop}.
\end{example}

For any bounded geometric antistructure,
 notice that
the action of $\theta_M$ on $M$ takes $H${--}orbits
to $H${--}orbits since $\theta_M(g\cdot m) = \theta_H(g)\cdot
\theta_M(m)$.
Let $M_{(H,\theta)}$ denote the subset consisting of $H${--}orbits
in $M$ which are \emph{fixed} by the $\theta${--}action.
Then  $M_{(H,\theta)} = \{ m\in M \Vertical \theta_M(m) \in H\cdot m\}$.
Note that $M_{(H,\theta)}$  is a $H${--}invariant subspace of $M$.
\begin{theorem}\label{optrick}
Suppose that $(\ccat{M,H}{R},\alpha,u)$ has a bounded geometric
antistructure $(\alpha,u)$
given by $(\theta_H, \theta_M, b, w)$, such that:
\begin{enumerate}
\renewcommand{\labelenumi}{(\roman{enumi})}
\item $M = O(K)$, where  $K$ is a finite $H${--}CW complex
and $M$ has the cone of the given $H${--}action on $K$,
\item $\theta_M$ is induced by a simplicial map on $K$,
\item $M_{(H,\theta)} \subseteq O(L):=N$ for some
$H${--}invariant subcomplex $L \subset K$, and
\item for some $i \leqq 1$, $I \subset K_i(\bCat{M}{H}{N}(R))$
 is a subgroup with $H^*(I)=0$.
\end{enumerate}
Then
$L_n^I(\bCat{M}{H}{N}(R),\alpha,u)=0$
for all $n$.
\end{theorem}
\begin{corollary}
Let $G=\Cy{2^rq}$, $q$ odd, be a cyclic group and
 $H \subset G$ the subgroup of index 2.
 Let $W$ be
a $G${--}representation, and
$N=\bigcup \{W^K \mid [G:K]\ \text{ is\ odd}\}$.
Then $L_n^I(\bCat{W}{H}{N}(\bz),\alpha,u)(q)=0$
on the  top component, where $(\alpha, u)$ is
the antistructure given above.
\end{corollary}
\begin{proof}[The proof of Theorem \ref{optrick}]
We extend the given $H${--}action on $M$ to a
simplicial action
of $G=\langle H,t\mid t^{-1}ht=\theta_H(h), t^2=b\rangle$ as described
above.
The proof is by induction on cells,
so suppose that $K$ is obtained
from $L$ by attaching exactly one $G${--}equivariant $k${--}cell
$D^k \times G/G_0$.
Since $M_{(H,\theta)}\subset O(L)=N$, it follows that
$G_0 \subset H$ and we may write $G/G_0 = H/G_0\,\sqcup \, tH/G_0$.
Now we define $M_1=O(L \cup (D^k\times H/G_0))$ and consider
the category $\ca = \Ccat{M_1,H}{O(L)}{R}$.
By construction, we have
\[\Ccat{M,H}{N}{R} = \ca \times  \ca^{op}\]
which has trivial $L${--}theory by Lemma \ref{Rop}.
Since the Tate cohomology of the $K_1$
decoration $I$ vanishes, we get
$L_n^I(\bCat{M}{H}{N}(R),\alpha,u)=0$.
\end{proof}

\section{ Calculations in bounded $K$-theory}
\label{boundedK}
We begin to compute the bounded transfers $\trf{W}$ by considering
the bounded $K${--}theory analogue.
In this section, $G = \Cy{2^rq}$ is cyclic of
order $2^rq$, with $r\geqq 2$ and $q\geqq 1$ odd.
By  \pone, Theorem \ref{twopoweriso} we can restrict our attention to those $W$
where the isotropy subgroups have $2${--}power index.
Let $G_i\subseteq G$ denote the subgroup
of index $[G:G_i]=2^{i}$ for $i=0, 1,\dots, r$.
As above, we reserve the notation $H < G$ for
the subgroup of index 2.

Any real, orthogonal $G${--}representation $W$
can be decomposed uniquely into isotypical
direct summands indexed by the subgroups
$K\subseteq G$, where in each summand $G$
operates  with isotropy group $K$ away from
the origin.
Since we assume that $W$ has isotropy of $2${--}power
index, we can write
\[W = W[0]\oplus W[1]\oplus \dots \oplus W[r]\]
where $W[i]$ is isotypic with isotropy group
$G_i$.
Thus $W[0]=\br^k$ is a trivial $G${--}representation,
and $W[1]$ is a sum of $\rminus$ factors.
We say that $W$ is \emph{complex} if $\dim W[0]$ and 
$\dim W[1]$ are even (in this case, $W$ is the underlying real
representation of a complex representation).
If $W$ is complex, then
 $\Wmax \subseteq W$  denotes a 
complex sub-representation of real $\dim\leq 2$
with maximal proper isotropy subgroup. If $W=W[0]$ then 
$\Wmax=0$.
Then $\Wmax$ is either irreducible
or $\Wmax = \rminus\oplus\rminus$.

We study bounded $K${--}theory by means of equivariant
filtrations of the control space.
The basic sequence is (see \cite{hp1}):
\[ \dots \to K_{i+1}(\brelcat{V}{U}) \RA{\bd_{i+1}} K_i(\bcat{U}) \to
K_i(\bcat{V}) \to K_i(\brelcat{V}{U}) \to \dots\]
valid for $U\subseteq V$ a closed $G${--}invariant subspace.
If $W_1$ is a complex representation with
 $\dim W_1=2$ and   isotropy group $K\neq G$,
let $U =\Uell$ be the union of $[G:K]$ rays from the origin
in $W_1$, which are freely permuted by $G/K$.
Then $W_1\setminus \Uell$ is a disjoint union of open fundamental
domains for the free $G/K${--}action. 
If $W = W_1\oplus W_2$, we call
$W_2 \subset \Uell\times W_2 \subset W$
the \emph{orbit type} filtration of $W$.

Recall that $t$ denotes a generator of $G$, and thus acts as an
isometry on the control spaces $M$ we use in the bounded categories
$\bcat{M}$.
Let $t_{\ast}$ denote the  action
 of $t $ on bounded $K${--}theory  induced by
 its action on the control space.

\begin{lemma}\label{orb}
Let $W = W_1\oplus W_2$, where $W_1$ is
a complex $2${--}dimensional
 sub-representation of $W$ with
minimal isotropy subgroup $ K \neq G$.
Then
 \begin{align*}
&K_{i+1}(\brelcat{W}{\Uell\times W_2})
\cong K_{i-1-k}(\ZK)\\
&K_{i}(\brelcat{\Uell\times W_2}{W_2})
\cong K_{i-1-k}(\ZK),
\end{align*}
where $k=\dim W_2$.
The boundary map $\bd_{i+1} = 1 - t_{\ast}$ in
the long exact sequence
of the orbit type filtration for $W$.
\end{lemma}
\begin{proof}
The bounded category $\brelcat{\Uell\times W_2}{W_2}$
of germs away from $W_2$  has effective fundamental
group $K$, as defined in \cite[3.13]{hp1}.
It therefore has the same $K${--}theory as
$\cc_{\br^{k+1}}(\cc_{pt}(\ZK))$.
The other case is similar.
The identification of $\bd_{i+1}$ with $1-t_{\ast}$ is
discussed in detail in the proof of Proposition \ref{boundary}.
\end{proof}
Since $\Kminus{j}(\ZK) = 0$ for $j\geqq 2$ by \cite{carter1},
this Lemma gives vanishing results for bounded $K${--}theory
as well.
\begin{lemma}\label{five}
Suppose that $W$ is complex, and $W^G=0$.
Then the inclusion map induces an isomorphism
\[K_i(\cc_{\Wmax\times \rminus,G}(\bz)) \to
K_i(\cc_{W\times \rminus,G}(\bz))\]
for $i \leqq 1$.
\end{lemma}
\begin{proof}
This is an argument using the orbit type filtration.
Let $W = \Wmax\oplus W_2$, and suppose that $W_2\neq 0$
or equivalently $\dim W_2\geqq 2$, since $W$ is complex.
Write $W = W' \oplus W''$ with $W'\subseteq W_2$,
 $\dim W' = 2$ and $\Iso(W') = K$ minimal.
We  choose $W''$ containing $\Wmax$,
and by induction we assume the result holds for $W''$.

Then applying the first part of Lemma \ref{orb}, we get the calculations
\eqncount
\begin{equation}K_i(\brelcatrminus{\Uell\times W''}{W''})
=K_{i-2-|W''|}(\ZK)\label{filta}
\end{equation}
and
\eqncount
\begin{equation}K_i(\brelcatrminus{W}{\Uell\times W''})
=K_{i-3-|W''|}(\ZK)\ .\label{filtb}
\end{equation}
Since  $\dim W'' \geqq 2$, we get the vanishing results
$K_i(\brelcatrminus{\Uell\times W''}{W''}) = 0$
for $i \leqq 2$ by \cite{carter1},
and
$K_i(\brelcatrminus{W}{\Uell\times W''}) =0$ for  $i\leqq 3$.
From the filtration sequence, it follows that
$K_i(\brelcat{W\times\rminus}{W''\times \rminus}) =0$
for $i\leqq 2$, and therefore
\[K_i(\bcat{W''\times\rminus})
\RA{\approx} K_i(\bcat{W\times\rminus})\]
for $i\leqq 1$.
We are done, by induction.
\end{proof}
\begin{corollary}\label{maxk}
$K_i(\brelcat{W\times\rminus}{\Wmax\times\rminus}) =0$
for $i \leqq 2$.
\end{corollary}
\begin{proof}
We continue the notation from above, and look at part of the
filtration sequence
\[  K_i(\brelcat{W''\times\rminus}{\Wmax\times\rminus})
 \to K_i(\brelcat{W\times\rminus}{\Wmax\times\rminus})
\to  K_i(\brelcat{W\times\rminus}{\Wmax\times\rminus})\ .\]
The first term is zero for $i\leqq 2$  by induction on dimension, and the
third term is zero for $i\leqq 2$ as above.
\end{proof}
We can obtain a little sharper result with some additional work.
First a useful observation:
\begin{lemma}\label{flipsign} Let $\ca$ be an additive category (with involution).
Then the  map $\br\to\br$ sending $x$ to $-x$ induces minus the identity on
$K$-theory (and $L$-theory) of $\cc_{\br}(\ca)$.
\end{lemma}
\begin{proof}
The category $\cc_{\br}(\ca)$ is filtered by the full subcategory
whose objects have support in a bounded neighborhood of $0$.
This subcategory is
equivalent to $\ca$ and the quotient category may be identified with
$\Ccat{[0,\infty)}{0}{\ca}\times \Ccat{(-\infty,0]}{0}{\ca}$
via the projection maps in an obvious way.

Consider the diagram
\[
\xymatrix{\ca\ar@{=}[d]\ar[r]&\cc_{\br}(\ca)\ar[d]\ar[r]&
\cc_{[0,\infty)}^{>0}(\ca)\times \cc _{(-\infty,0]}^{>0}(\ca)\ar[d]\\
\ca\ar[r]&\cc_{[0,\infty)}(\ca)\ar[r]&\cc_{[0,\infty)}^{>0}
}
\]
where the vertical map is induced by $x\mapsto |x|$.
 In the lower horizontal row, $K$ and $L${--}theory of the
 middle term is trivial, so the
boundary map will be an isomorphism.
The lower row splits off the upper row in
two different ways, one induced by including $[0,\infty)\subset \br$ and the
other by sending $x\in[0,\infty)$ to $-x\in\br$.
Under these two splittings we
may identify $K$ or $L${--}theory of the quotient
$\cc_{[0,\infty)}^{>0}(\ca)\times \cc_{(-\infty,0]}^{>0}(\ca)$ with
$\cc_{[0,\infty)}^{>0}\times \cc_{[0,\infty)}^{>0}$ and under this
identification, the flip map of $\br$ corresponds to interchanging the two
factors.
On $K${--}theory (or $L${--}theory)
we conclude that the exact sequence is
of the form
\[ 0\to A_*\to A_*\times A_*\RA{+} A_*\to 0\ .\]
The flip action on the last term is trivial, and on the middle
term it interchanges the two factors, so the inclusion must send $a$ to
$(a,-a)$.
Hence the flip action on the first term must be $a\mapsto -a$.
\end{proof}

As above, $t_{\ast}$ denotes the induced action
 of a generator $t \in G$
on $K${--}theory, and $\varepsilon\colon G \to \orient$ the non-trivial  action
of $G$ on $\rminus$.

\begin{proposition}\label{boundary}
Let $W$ be a
complex $2$-dimensional $G${--}representation with $W^G=0$,
and isotropy subgroup $K$.
Then under the isomorphisms of (\ref{filta}) and (\ref{filtb})
the complex
\[ K_{i+1}(\brelcatrminus{W}{\Uell})
\RA{\bd_a} K_{i}(\brelCatrminus{\Uell}) \RA{\bd_b}
K_{i-1}(\brelcat{\rminus}{0})\]
with $\bd_b\circ\bd_a =0$ is isomorphic to
\[K_{i-2}(\ZK) \RA{\bd'_a} K_{i-2}(\ZK) \RA{\bd'_b} K_{i-2}(\ZH)\]
where
$\bd'_a = 1-t_{\ast} =2$ and
$\bd'_b = \Ind_H\circ (1-t_{\ast}) = 0$.
\end{proposition}
\begin{proof}
The orbit type filtration is based on a $G${--}equivariant
simplicial model for $W$, where $G$ acts through the
projection to $G/K$.
The third term in the complex is $K_{i-1}(\brelcat{\rminus}{0})
=K_{i-2}(\ZH)$ and the identification of the  boundary maps
follows from the definition of the germ categories.

To compute $\bd_a$, we use the isomorphisms between the domain of
$\bd_a$  and $K_{i-2}(\bz K)$  obtained by noticing that every element is induced
from an  element  of
\[K_{i+1}(\Ccat{C\times \br,K}{\bd C\times \br}{\bz)},\]
where $C$ is the region between two adjacent half{--}lines of $\Uell$.
This follows since the regions in the complement of $\Uell\times\rminus$
 are disjoint and the boundedness condition ensures there is no
interference.
 Similarly the isomorphism of the range of $\bd_a$  with
$K_{i-2}(\bz)$ is obtained by noticing that every element is induced from
$\cc_{h\times \br,K}^{>\br}$ where $h$
 is just one half{--}line in $\Uell$, and we can think of
$\bd C=h\cup th$.
To compute the boundary we first take the standard
boundary to $\cc_{\bd C\times \br,K}(\bz)$ which is an isomorphism, and then
map away from $0\times \br$.
 It follows from the proof of Lemma \ref{flipsign} above
that this map is of the form $a \mapsto (a,-a)$ in $K$ or $L${--}theory.
In this picture, the support of one of the boundary components
is along $h$ and the other along $th$.
We need to use the group action to associate both elements to
the same ray.
Since $t$ flips the $\br${--}factor, we get a change of sign before
adding, so $t_\ast = \varepsilon(t) = -1$ and  $\bd_a$ sends $a$ to $2a$.

To compute $\bd_b$, we start with an element in the source of $\bd_b$ which
as above is identified via induction with
$K_i(\brelcat{h\times \br,K}{\br})$.
The boundary first sends this isomorphically to $\cc_{\br,K}(\bz)$,
 then by induction to $\bCatrminus$,  and
then via the natural map to the range of $\bd_b$, which is
$\cc_{\rminus,G}^{>0}$.
We have a commutative square
\[
\xymatrix{K_1(\Ccat{h \times \br,K}{\br}{\bz})
\ar[r]^{\Ind} \ar[d]_{\bd}
&K_1(\brelCatrminus{\Uell})\ar[d]^{\bd_b}\\
K_0(\Ccat{\br,K}{0}{\bz}) \ar[r]^{\Ind}&K_0(\Ccat{\rminus,G}{0}{\bz})
}
\]
where $\Ccat{\br,K}{0}{\bz} =\cc_{\br,K}(\bz)\times
\cc_{\br,K}(\bz)$.
Under this identification, the natural map to the germ category
$K_0(\cc_{\br,K}(\bz)) \to K_0(\Ccat{\br,K}{0}{\bz})$
is just $a \mapsto (a,-a)$, and
the induction map
\[\Ind\colon K_0(\Ccat{\br,K}{0}{\bz}) \to K_0(\Ccat{\rminus,G}{0}{\bz})\]
is  given by $(a,b)\mapsto \Ind_H(a+t_{\ast} b)$.
In this case the action of $t$ on
$\cc_{\rminus,G}(\bz)$ is  the identity since any element is
invariant under the action of $G$, hence under the action of $t$.
It follows that $\Ind(a,-a) = \Ind_H(a-t_{\ast}a) =0$ as required.
\end{proof}

\begin{lemma}\label{boundtwo}
Let $W$ be a complex $2$-dimensional $G${--}representation
with proper isotropy group $K$.
The boundary map
\[K_{i+1}(\brelCatrminus{\Uell}) \to K_i(\bCatrminus)\]
is zero for $i\leqq 1$.
\end{lemma}
\begin{proof}
If $i \leqq -1$ the domain of this boundary map is zero, so the result is trivial.
For $i=1$ we use the injection $\wh(\bCatrminus) \to \wh(\ccat{\rminus,G}{\zadictwo})$,
which follows from the vanishing of $SK_1(\ZG)$. But
$\wh(\ccat{\rminus,G}{\zadictwo}) = \wh( \Zadictwo{G})/\wh(\Zadictwo{H})$ and 
the group $ \wh( \Zadictwo{G}) = \prod_{d\mid q} \zadictwo [\zeta_d] G_2$,
where $G_2 \subset G$ is the $2$-Sylow subgroup and $q$ is the odd part of the order
of $G$. Now Oliver \cite[Thm. 6.6]{o1}  constructed a short exact sequence
$$1 \to \wh(\zadictwo [\zeta_d] G_2) \to 
\zadictwo [\zeta_d] G_2 \to \langle -1\rangle \times G_2 \to 1$$
by means of the integral $2$-adic logarithm. This sequence is natural with respect to
inclusion of subgroups, so we may use it to compare $ \wh( \Zadictwo{G})$ and
$ \wh( \Zadictwo{H})$. Since each corresponding term injects, and the middle quotient is 
$\bz$-torsion free, we conclude that $\wh(\ccat{\rminus,G}{\zadictwo})$ is
 also $\bz$-torsion free. Since $K_{2}(\brelCatrminus{\Uell} = K_{0}(\ZK)$ is
torsion (except for $\bz = K_0(\bz)$ which is detected by projection to the
trivial group),
the given boundary map is zero.

For $i=0$, we use the surjection $K_0(\qlocal K) \to \Kminus{1}(\ZK)$,
and compute  with $\qlocal$
coefficients and $i=1$.
We will list the steps, and leave the details to the reader.
First, compute that $K_1(\ccat{\rminus,G}{\qlocal})=K_1(\qlocal G)/
K_1(\qlocal H)$ surjects onto $K_1(\ccat{\Uell\times\rminus,G}{\qlocal})$,
by means of a braid containing the cone point inclusions
into $\ccat{\rminus,G}{\qlocal}$ and $\ccat{\Uell\times\rminus,G}{\qlocal}$.
Second, prove that $K_1(\ccat{\Uell\times\rminus,G}{\qlocal})$
fits into a short exact sequence
\[0 \to K_1(\ccat{\Uell,H}{\qlocal}) \to K_1(\ccat{\Uell,G}{\qlocal})
\to  K_1(\ccat{\Uell\times\rminus,G}{\qlocal}) \to 0\]
by means of a braid containing the inclusion
$\Uell \subseteq \Uell\times\rminus$.
Finally, compute the first two terms
$K_1(\ccat{\Uell,G}{\qlocal})=K_1(\qlocal G)/K_1(\qlocal K)$,
and $K_1(\ccat{\Uell,H}{\qlocal})=K_1(\qlocal H)/K_1(\qlocal K)$
by comparing the groups under the inclusion $H < G$.
We conclude that
\[K_1(\ccat{\rminus,G}{\qlocal}) \cong
K_1(\ccat{\Uell\times \rminus,G}{\qlocal})\] under the inclusion map,
and hence $\bd =0$.
\end{proof}

\begin{corollary}\label{kersimple}
Let $W$ be a complex $G${--}representation with $W^G=0$.
Then the inclusion induces  an isomorphism
$K_i(\bCatrminus) \to K_i(\bcatrminus{W})$
for $i\leqq 0$, and an injection for $i=1$.
If $\Kminus{1}(\ZK)=0$ for the maximal proper isotropy
 groups $K$ of $W$,
then $K_1(\bCatrminus) \cong K_1(\bcatrminus{W})$.
\end{corollary}
\begin{proof}
We may assume that $\dim W=2$, and apply the filtering
argument again.
By (\ref{filta}) and (\ref{filtb}) we get
$K_i(\brelCatrminus{W}) =0$ for $i\leqq 0$, and
$K_1(\brelCatrminus{W})$ is a quotient of
$K_1(\brelCatrminus{\Uell})=\Kminus{1}(\ZK)$.
Since the composition
\[K_{1}(\brelCatrminus{\Uell})
\to K_{1}(\brelCatrminus{W}) \to K_0(\bCatrminus)\]
is zero  by Lemma \ref{boundtwo}, the result follows for $i\leqq 0$.

Similarly, $K_2(\brelCatrminus{W})$
is a quotient of $K_2(\brelCatrminus{\Uell}) = K_0(\ZK)$,
because the boundary map $K_2(\brelcatrminus{W}{\Uell})=\Kminus{1}(\ZK)$
to $K_1(\brelCatrminus{\Uell})=\Kminus{1}(\ZK)$ is multiplication
by $2$, and hence injective.
Then we make the same argument, using Lemma \ref{boundtwo}. 

If we also assume $\Kminus{1}(\ZK)=0$, then
$ K_1(\brelCatrminus{W}) =0$ so we get the isomorphism
$K_1(\bCatrminus) \RA{\approx} K_1(\bcatrminus{W})$.
\end{proof}

\section{ The double coboundary}
\label{doublecob}
The composite $\Delta$
of maps from the $L^h-L^p$ and $L^s -L^h$ Rothenberg
sequences
\[
\xymatrix{
H^{i}(\widetilde K_0(\ZG))\ar[r]\ar[dr]^{\Delta} & \Lgroup{h}{i}{\ZG}\ar[d]\\
&  H^i(\wh(\ZG))}
\]
(see \cite{ht1}) has an algebraic description
by means of a ``double coboundary" homomorphism
\[\delta^2\colon H^{i}(\widetilde K_0(\ZG))\to H^i(\wh(\ZG))\]
In this section, we will give a brief description due to
Ranicki \cite{ra13} of this homomorphism
(see also \cite[\S 6.2]{ra8} for related material on
``interlocking" exact sequences in $K$ and $L${--}theory).

Let $X$ be a space with a $\cy 2$ action $T\colon X \to X$,
and define homomorphisms
\[\Delta\colon H^i(\pi_n(X)) \to H^i(\pi_{n+1}(X))\]
by sending $g\colon S^n \to X$ to
\[h\cup (-1)^iTh\colon S^{n+1} = D^{n+1}_{+} \cup_{S^n} D^{n+1}_{-}
\to X\]
for any null-homotopy $h\colon D^{n+1} \to X$ of the map
$g + (-1)^{i+1}Tg\colon S^n \to X$.

The maps $\Delta$ lead to a universal description of double
coboundary maps, as follows.
Let $f\colon X \to Y$ be a $\cy 2${--}equivariant map of
spaces with $\cy 2$ action, and consider the long exact
sequence
\[\dots \to \pi_n(X) \RA{f} \pi_n(Y) \to \pi_n(f) \to \pi_{n-1}(X)
\to \pi_{n-1}(Y) \to \dots\]
We define
$I_n=\ker(f\colon \pi_n(X) \to \pi_n(Y))$ and
$J_n=\Image(\pi_n(Y) \to \pi_{n}(f))$, and get an exact
sequence
\[0 \to \pi_n(X)/I_n \to \pi_n(Y) \to \pi_n(f) \to I_{n-1} \to 0\]
which can be spliced together from the short exact sequences
\eqncount
\begin{equation}
\begin{split} \label{tateseq}
&0 \to  \pi_n(X)/I_n \to \pi_n(Y) \to J_n \to 0\\
&0 \to J_n \to \pi_n(f) \to I_{n-1} \to 0\ .
\end{split}
\end{equation}

Then it follows directly from the definitions that the double
coboundary
\[\delta^2\colon H^i(I_{n-1}) \RA{\delta} H^{i+1}(J_n) \RA{\delta}
H^i(\pi_n(X)/I_n)\]
from the Tate cohomology sequences induced
by (\ref{tateseq}) is given by the composite
\[\delta^2\colon H^i(I_{n-1}) \RA{inc_{\ast}} H^i(\pi_{n-1}(X))
\RA{\Delta} H^i(\pi_n(X)) \RA{proj_{\ast}} H^i(\pi_n(X)/I_n)\ .\]
If we can pick the map $f\colon X \to Y$ appropriately,
say with $I_n=0$ and $I_{n-1} = \pi_{n-1}(X)$, this gives
an algebraic description of $\Delta$.

In later sections we will use the relative Tate cohomology groups
$H^i(\Delta)$, which are just (by definition) the relative Tate
cohomology groups \cite[p.166]{ra8} of the map $\pi_n(Y) \to \pi_n(f)$ in
the long exact sequence above.
These groups fit into the commutative braids given in \cite{ra13}
which will be used in the proofs of Theorems A-C.

We now give some examples, with $G$ denoting a finite cyclic
group as usual. These arise from homotopy groups of certain fibrations
of algebraic $K$-theory spectra.
\begin{example}
There is an exact sequence \cite{ht1}
\[0 \to \wh(\ZG) \to \wh(\Zadictwo{G}) \to  \wh(\Zreltwo{G})
\to \widetilde K_0(\ZG) \to 0\]
of $\cy 2$ modules and the associated double coboundary in Tate
cohomology equals
\[\Delta\colon H^{i}(\widetilde K_0(\ZG)) \to H^i(\wh(\ZG))\ .\]
The point here is that $\ker(\wh(\ZG) \to \wh(\Zadictwo{G}))=0$
\cite{o1}, and the map $\widetilde K_0(\ZG) \to
\widetilde K_0(\Zadictwo{G})$ is zero by a result of Swan
\cite{sw1}.
We could also use the exact sequence
\[0 \to  \wh(\ZG) \to \wh(\zlocal G)\oplus \wh(\bq G)\to
 \wh(\qlocal G)) \to \widetilde K_0(\ZG) \to 0\]
to compute the same map $\Delta$.
\end{example}
\begin{example}
There is an exact sequence
\begin{multline*}
0 \to \wh(\bCatrminus) \to \wh(\ccat{\rminus, G}{\zlocal\oplus \bq})
\to \\
\wh(\ccat{\rminus,G}{\qlocal})  \to \widetilde K_0(\bCatrminus) \to 0
\end{multline*}
where $\wh(\bCatrminus) = \wh(\rel{\ZH}{\ZG})$ and
$\widetilde K_0(\bCatrminus) = \widetilde K_0(\rel{\ZH}{\ZG})$.
The injectivity on the left follows because
 $K_2(\ccat{\rminus,G}{\qlocal})$ is a quotient of $K_2(\qlocal G)$,
mapping trivially through $K_1(\ZG)$ into $\wh(\bCatrminus)$
(since $SK_1(\ZG)=0$ \cite{o1}).
We therefore get an algebraic description of $\delta^2\colon
H^i(\widetilde K_0(\rel{\ZH}{\ZG^{-}})) \to H^i(\wh(\rel{\ZH}{\ZG^{-}}))$
as used in the statement of Theorem B.
\end{example}

\begin{example}
There is an exact sequence
\begin{multline}
0 \to \wh(\bcatrminus{W}) \to
 \wh(\ccat{W\times \rminus,G}{\zlocal \oplus \bq})
\to\\
 \wh(\ccat{W\times \rminus,G}{\qlocal})
 \to \widetilde K_0(\bcatrminus{W}) \to 0
\end{multline}
for any complex $G${--}representation $W$ with $W^G=0$.
We therefore get an algebraic description of $\delta^2\colon
H^i(\widetilde K_0(\bcatrminus{W})) \to H^i(\wh(\bcatrminus{W}))$
as used in the statement of Theorem C.
\end{example}
\begin{lemma}
For complex
 $G${--}representations $W_1\subseteq W$ with $W^G=0$,
there is a commutative diagram
\[
\xymatrix{H^i(\widetilde K_0(\bcatrminus{W_1}))
\ar[d]_{c_{\ast}}\ar[r]^{\delta^2}
& H^i(\wh(\bcatrminus{W_1}))\ar[d]^{c_{\ast}}\\
H^i(\widetilde K_0(\bcatrminus{W})) \ar[r]^{\delta^2}
& H^i(\wh(\bcatrminus{W}))
}
\]
where the vertical maps are induced by the inclusion $W_1\subseteq W$.
\end{lemma}
For our applications, the main point of the double coboundary
description is that it permits these maps induced by
cone point inclusions to be computed using bounded $K${--}theory,
instead of bounded $L${--}theory.

The double coboundary maps also commute with restriction
to  subgroups of $G$.
\begin{proposition}\label{ninetwo}
Let $G_1<G$ be a subgroup of odd index, and $H_1<G_1$ have index $2$,
then there are twisted restriction maps
$$H^i(\widetilde K_0(\rel{\ZH}{\ZG^-})) \xrightarrow{\Res}
 H^i(\widetilde K_0(\rel{\ZH_1}{\ZG_1^-}))
$$
 and
 $$H^i(\wh(\rel{\ZH}{\ZG^-})) \xrightarrow{\Res}
H^i(\wh(\rel{\ZH_1}{\ZG_1^-}))$$
such that the diagram
\[
\xymatrix{H^i(\widetilde K_0(\rel{\ZH}{\ZG^-}))
\ar[r]^{\delta^2}\ar[d]_{\Res}&
H^i(\wh(\rel{\ZH}{\ZG^-}))\ar[d]^{\Res}\\
H^i(\widetilde K_0(\rel{\ZH_1}{\ZG_1^-}))\ar[r]^{\delta_1^2}&
  H^i(\wh(\rel{\ZH_1}{\ZG_1^-}))
}
\]
commutes.
\end{proposition}
\begin{proof}
The vertical maps are twisted restriction maps given by
composing
the  twisting isomorphisms 
$H^i(\widetilde K_0(\rel{\ZH}{\ZG^-})) \cong
 H^{i+1}(\widetilde K_0(\rel{\ZH}{\ZG}))$ 
and
$H^i(\wh(\rel{\ZH}{\ZG^-}) \cong H^{i+1}(\wh(\rel{\ZH}{\ZG})$,
discussed in \pone, Section \ref{rmin},
with the  restriction maps induced by the inclusion
$\ZG_1\to \ZG$ of rings with involution.
Since $G_1 <G$ has odd index, $H_1 = H\cap G_1$ and the composition
$\Res_{G_1}\circ \Ind_{H}$ lands in the image of $\Ind_{H_1}$
by the double coset formula.
\end{proof}
This can be generalized to the double coboundary maps
used in the statement of Theorem C, under certain
conditions.
\begin{proposition}\label{ninethree}
Let $G_1<G$ be an odd index subgroup, and $H_1<G_1$ have index $2$.
Suppose that
$W$ only has proper isotropy subgroups  of $2${--}power index.
Then there are a twisted restriction maps
$H^i(\widetilde K_0(\bcatrminus{W}))\to
H^i(\widetilde K_0(\ccat{\Res W\times\rminus,G_1}{\bz}))$
and
$H^i(\wh(\bcatrminus{W})) \to
 H^i(\wh(\ccat{\Res W\times\rminus,G_1}{\bz}))$
which commute with the corresponding
 double coboundary maps $\delta^2_W$ and $\delta^2_{\Res W}$.
\end{proposition}

\section{ Calculations in bounded $L${--}theory}
 \label{boundedL}

Suppose that $\sigma=\sigma(f) \in L^h_0(\ZG)$
is the surgery obstruction arising from a normal cobordism
between $S(V_1)$ and $S(V_2)$, as in the statement of
\pone,  Theorem \ref{maintest}.
In this section, we establish two important properties of
$\trf{W}(\sigma)$ in preparation for the proof of Theorem C.
Unless otherwise mentioned, all bounded categories will
have the standard orientation (see Example \ref{standard}).

For a complex $G${--}representation $W$ the standard
orientation is trivial, and
the cone point inclusion  $0 \in W$ induces the map
\[c_*\colon  L_n^h(\ZG) \to L_{n}^h(\bcat{W})\]
Note that the presence
of an $\rminus$ factor introduces a non-trivial orientation
at the cone point
\[c_*\colon  L_{n+1}^h(\ZG, w) \to  L_{n+1}^h(\bcatrminus{W})\]
where $w\colon G \to \orient$ is the non-trivial
projection.
The properties are:
\begin{theorem}\label{one}
Suppose  $W$ is a complex $G${--}representation with
no  $\rplus$ summands.
If $\trf{W\times \rminus}(\sigma)
\in L^h_{2k+1}(\bcatrminus{W})$
is a torsion element, 
then
\begin{enumerate}
\renewcommand{\labelenumi}{(\roman{enumi})}
\item there exists a torsion element $\hat\sigma\in 
L_{2k+1}^h(\bCatrminus)$ such that 
\[c_*(\hat\sigma) = \trf{W\times \rminus}(\sigma)
\in  L^h_{2k+1}(\bcatrminus{W})\]
\item there exists a torsion
 element $\hat\sigma \in L^p_{2k+1}(\ZG, w)$
such that
\[c_*(\hat\sigma) = \trf{W\times \rminus}(\sigma)
\in  L^p_{2k+1}(\bcatrminus{W})\]
where $\dim W = 2k$.
\end{enumerate}
\end{theorem}
We remark that the condition ``$\trf{W\times \rminus}(\sigma)$
is a torsion element"  follows from  the assumption
$\Res_H (V_1\oplus W) \oplus \rplus
 \sim_t \Res_H (V_2\oplus W)\oplus \rplus$ in 
Theorem C,
as an immediate consequence of Proposition \ref{trftorsion}.
Before giving the proof,  we need some preliminary
results.
When the $L${--}theory decoration is not explicitly given,
we mean $\Linfty$.

\begin{lemma}
Let $W$ be a complex $G${--}representation.
Then
\[ L_{2k+1}(\bcat{W})\otimes \bq = 0\]  for $k \geqq 0$
\end{lemma}
\begin{proof}
We argue by induction on $\dim W$, starting with
\[L_{2k+1}(\ZG) = \cy 2\oplus H^1(\Kminus{1}(\ZG))\]
 which is all $2${--}torsion.
It is enough to prove the result for the top component
$ L_{2k+1}(\bcat{W})(q)$, and therefore by \pone, 
Theorem \ref{twopoweriso} we may
assume that the isotropy groups of $W$ all have
$2${--}power index.
Since we are working with $\Linfty$, we may ignore
$\rplus$ summands of $W$.

Let $W = W'\oplus W''$ where $\dim W'=2$ and
$W'$  has minimal isotropy group $K$.
 We assume the result for $W''$,
and let $\Uell \subset W'$ be a $G$-invariant
set of rays from the origin, dividing $W'$ into fundamental
domains for the free $G/K$-action.

Then
\[L_{n}(\bCat{\Uell \times W''}{G}{W''}(\bz))
=L_{n-1}(\bcat{W''})
=L_{n-1-|W''|}(\ZK),\]
which is torsion for $n$ even, and
\[L_{n}(\bCat{W}{G}{\Uell \times W''}(\bz))
=L_{n-2}(\bcat{W''})
=L_{n-2-|W''|}(\ZK),
\]
which is torsion for $n$ odd.
Moreover, we have a long exact sequence
\[
\dots\to L_{n}(\bCat{\Uell \times W''}{G}{W''}(\bz)) \to
L_{n}(\bCat{W}{G}{W''}(\bz)) \to
L_{n}(\bCat{W}{G}{\Uell \times W''}(\bz)) \to \dots
\]
We claim that  the first map in this sequence is
rationally  injective.
The previous map in the long exact sequence
is
\[ L_{n+1}(\bCat{W}{G}{\Uell \times W''}(\bz))
\longrightarrow{\bd}
 L_{n}(\bCat{\Uell \times W''}{G}{W''}(\bz)),\]
which may be identified (using the isomorphisms above
and Proposition \ref{boundary})
with a geometrically induced map
\[L_{n-1-|W''|}(\ZK) \RA{1-u} L_{n-1-|W''|}(\ZK)\]
``multiplication by $1-w(t)$", where $\bar t = w(t) t^{-1}$
from the action of $t$ in the antistructure used to define the
$L$-group.

In the oriented case, $w(t)= +1$ and
this boundary map is zero.
Therefore, we conclude that
\[\xymatrix{L_{2k+1}(\bCat{\Uell \times W''}{G}{W''}(\bz))
 \ar[r]^{\approx}&
L_{2k+1}(\bCat{W}{G}{W''}(\bz))} \]
is a rational isomorphism
and so
\[L_{2k+1}(\bCat{W}{G}{W''}(\bz)) = L_{2k-|W''|}(\ZK)\ .\]
Finally, we will substitute this computation and our
inductive assumption into the long exact sequence
\begin{multline*}
 \dots \to L_{2k+1}(\bcat{W''}) \to
 L_{2k+1}(\bcat{W}) \to\\
L_{2k+1}(\bCat{W}{G}{W''}(\bz)) \to
L_{2k}(\bcat{W''})  \to
\end{multline*}
and obtain an exact sequence
\[0 \to L_{2k+1}(\bcat{W}) \to
L_{2k+1}(\bCat{W}{G}{W''}(\bz)) \to L_{2k}(\bcat{W''}).\]
However, the second map in this sequence can be identified
with the inclusion map
\[
\Ind_K\colon  L_{2k}(\cc_{W'', K}(\bz))\to  L_{2k}(\bcat{W''})\]
and the composition $\Res _K \circ  \Ind_K$ is multiplication
by $[G:K]$, which is a rational isomorphism.
Therefore $\Ind_K$ is injective
and $ L_{2k+1}(\bcat{W}) =0$.
\end{proof}

Another computation we will need is
\begin{lemma}\label{fourtors} $L_{n}(\brelCatrminus{W})\otimes \bq = 0$
for W a complex G-representation.
\end{lemma}
\begin{proof}
We may assume that $W^G=0$  and argue by
 induction on the dimension of W.
We write $W = W' \oplus W''$,
where $\Iso(W' ) = K$ is minimal ($2${--}power index isotropy subgroups
may be assumed as usual).
We have 2 long exact sequences (all $L${--}groups
are  tensored  with $\bq$):

(i)  from the inclusion $W''\oplus \rminus \subset W \oplus \rminus$.
For short, let $A = \brelCatrminus{W}$,
$B = \brelCatrminus{W''}$, and then
$A/B = \brelcatrminus{W}{W''}$.
We need the
piece of the $L$-group sequence:
\eqncount
\begin{equation}
 ... \to L_{n}(B) \to  L_{n}(A) \to  L_{n}(A/B) ... \label{foura}
\end{equation}
and note that $L_{n}(B) = 0$ by induction.

(ii) we study $L_{n}(A/B)$ by looking at the usual rays
 $\Uell \subset W'$ which divide the 2-dim representation
$W'$ into G/K chambers.
 Let
    \[D_0 =\brelcatrminus{\Uell\times W''}{W''}\]
and
 \[   D_1 =\brelcatrminus{W}{\Uell\times W''}\]
Then we need the $L${--}group sequence
\eqncount
\begin{equation} \to L_{n}(D_0) \to  L_{n}(A/B) \to  L_{n}(D_1) \to L_{n-1}(D_0)
\label{fourb}
\end{equation}
But the groups $L_{n}(D_i)$ are the ones we have been computing
 by using the chamber structure.
 In particular,
\[L_{n}(D_0) = L_{n-2}(\cc_{W'',K}(\bz))\]
and
\[L_{n}(D_1) = L_{n-3}(\cc_{W'',K}(\bz)).\]
But since $K$ is the minimal isotropy group in $W$,
it acts trivially on $W''$
and these $L$-groups are just
\[
        L_{n}(D_0) = L_{n-2-|W''|}(\ZK) \]
and
\[L_{n}(D_1) = L_{n-3-|W''|}(\ZK).\]
Now the boundary map
\eqncount
\begin{equation}
        L_{n}(D_1) \to L_{n-1}(D_0)
\end{equation}
in the sequence is just multiplication by
 $1-w(t)=2$ ($w(t)=-1$ since $\rminus$
is non-oriented), and
this is a rational isomorphism.
Therefore, $L_{n}(A/B) = 0$ by (\ref{fourb}) and
substituting back into (\ref{foura}), the  $L_n(A)=0$.

\end{proof}
Now a more precise result in a special case:
\begin{lemma}\label{four}
Let W be a complex $2${--}dimensional G-representation
with $W^G=0$.
Then $L^s_{3}(\brelCatrminus{W}) = 0$
\end{lemma}
\begin{proof}
We have an exact sequence
\begin{multline*}
\dots \to L^s_3(\brelCatrminus{\Uell})
\to L^s_3(\brelCatrminus{W})\to\\
 L^I_3(\brelcatrminus{W}{\Uell})
\to L^s_2(\brelCatrminus{\Uell})
\to \dots
\end{multline*}

arising from the orbit type filtration.
Here
\begin{align*}
I &= \Image(K_2(\brelCatrminus{W}) \to K_2(\brelcatrminus{W}{\Uell}))\\
&=\ker(K_2(\brelcatrminus{W}{\Uell}) \to K_1(\brelCatrminus{\Uell}))\\
&=\ker(\Kminus{1}(\ZK) \RA{2} \Kminus{1}(\ZK))
\end{align*}
by Proposition \ref{boundary}, where $K$ is the isotropy subgroup of $W$.
But  $\Kminus{1}(\ZK)$ is torsion{--}free,  so $I = \{ 0 \} \subseteq K_2$.
Now substitute the computations
 \begin{align*} &L^s_3(\brelCatrminus{\Uell}) = L_1^p(\ZK) =0\\
&L^I_3(\brelcatrminus{W}{\Uell}) = L^p_0(\ZK)\\
&L^s_2(\brelCatrminus{\Uell}) = L^p_0(\ZK)
\end{align*}

into  the exact sequence.
The boundary map
\[L^I_3(\brelcatrminus{W}{\Uell})
\to  L^s_2(\brelCatrminus{\Uell})\]
 is multiplication by $2$ so
$ L^s_3(\brelCatrminus{W})=0$.
\end{proof}
\begin{remark}\label{fourex}
 The same method shows that
$ L^s_1(\brelCatrminus{W})=L_3^p(\ZK)=\cy 2$
for $\dim W=2$ as above, assuming that $K\neq 1$.
\end{remark}

Our final preliminary result is a Mayer{--}Vietoris sequence:
\begin{lemma}\label{two} Let $W$ be a complex $G${--}representation
with $W^G=0$.
Let $W = W_1 \oplus W_2$ be a direct
sum decomposition, where $W_1=\Wmax$.
\begin{enumerate}
\renewcommand{\labelenumi}{(\roman{enumi})}
\item
There is a long exact sequence
\begin{multline*}
\dots \to L^I_{n+1}(\brelcat{W\times \rminus}
{W_1\times \rminus\cup W_2\times \rminus})
\RA{\bd_{n+1}}
L^s_{n}(\brelCatrminus{W})
 \to \\
L^s_{n}(\brelcatrminus{W}{W_1})\oplus
 L^s_{n}(\brelcatrminus{W}{W_2})
 \to
L^I_{n}(\brelcat{W\times \rminus}
{W'\times \rminus\cup W''\times \rminus})
\end{multline*}

of bounded $L$-groups, where
\[I = \Image\left [ K_2(\brelcatrminus{W}{W_2}) \to
K_2(\brelcat{W\times \rminus}
{W'\times \rminus\cup W''\times \rminus})\right ]\]
is the decoration subgroup.
\item For $n\equiv 3\Mod{4}$, the boundary map
$\bd_{n+1} =0$.
\end{enumerate}
\end{lemma}
\begin{proof}
The Mayer{--}Vietoris sequence in $\Linfty$
\begin{multline*}
\dots \to L_{n+1}(\brelcat{M}{U\cup V}) \to
L_n(\brelcat{M}{U\cap V}) \to \\
L_n(\brelcat{M}{U})\oplus L_n(\brelcat{M}{V})\to
L_n(\brelcat{M}{U\cup V}
\end{multline*}
where $U$, $V$ are nice $G${--}invariant subspaces of
the control space $M$,
follows from an excision isomorphism
\[L_{n}(\cc^{>U}_{U\cup V,G}(\bz))
 = L_{n}(\cc^{>U\cap V}_{V,G}(\bz))\]
and standard diagram-chasing.
We apply this to $M = W \times \rminus$, $U = W_1\times \rminus$
and $V = W_2\times\rminus$, where $U\cap V = \rminus$.
The decorations follow from Section \ref{excision}
and Corollary \ref{maxk}.

To see that $\bd_{0}=0$, note that this boundary map
is the composition of
\[\bd\colon L^I_{0}(\brelcat{W\times \rminus}
{W'\times \rminus\cup W''\times \rminus}) \to
L^s_{3}(\brelcat{W'\times \rminus\cup W''\times \rminus}{W_2\times \rminus})\]
and an excision isomorphism
\[L^s_{3}(\brelcat{W'\times \rminus\cup W''\times \rminus}{W_2\times \rminus})
\cong L^s_3(\brelCatrminus{W_1})\ .\]
But $L^s_3(\brelCatrminus{W_1}) =0$ by Lemma \ref{four}.
\end{proof}

\begin{corollary}\label{twocor}
Suppose  $W$ is a  complex $G${--}representation with
no  $\rplus$ summands.
Then the image of
 $\trf{W\times \rminus}(\sigma)$
is zero in
$L^s_{2k+1}(\brelCatrminus{W})$,
where $\dim W =2k$.
\end{corollary}
\begin{proof}
First recall from \pone, Theorem \ref{ktrans} that
$\trf{W\times \rminus}(\sigma) \in L^{c_1}_{2k+1}(\bcatrminus{W})$
where $c_1=\Image(\wh(\ZG) \to \wh(\bcatrminus{W}))$
under the cone point inclusion.
Therefore, the image of $\trf{W\times \rminus}(\sigma)$
``away from $\rminus$" lands in $L^s_{2k+1}(\brelCatrminus{W})$.
We may therefore apply the Mayer{--}Vietoris sequence
from Lemma \ref{two}.

As usual, we  may assume that
$W$ has only $2${--}power isotropy subgroups,
and we will argue by induction on $\dim W$ and on the
orbit type filtration of $W$.
The result is true for $\dim W = 2$  (by Lemma \ref{four})
or for $W$ a free representation (by a direct calculation
following the method of Lemma \ref{fourtors}).
For sets of $2${--}power index
subgroups (ordered by inclusion, with repetitions allowed),
we say that $\{K'_1,\dots,K'_s\} < \{K_1,\dots, K_t\}$
if $s <t$ or $s=t$, $K'_i =K_i$ for $1\leqq i <j\leqq s$, and
$K'_j \subsetneqq K_j$.
Let $\Iso(W) = \{K_1,\dots,K_t\}$ denote the set of
isotropy subgroups of $W$ ordered by inclusion.
Define an ordering by $W' < W$ if (i) $\dim W' <\dim W$
or (ii) $\dim W'=\dim W$ but $\Iso(W') < \Iso(W)$.

Let $W = W_1\oplus W_2$, where $W_1$ has maximal
isotropy subgroup $K$, and suppose the result holds
for all $W' <W$.
If $\dim W\equiv 2\Mod{4}$, we get an injection
\[L^s_3(\brelCatrminus{W}) \to
 L^s_3(\brelcatrminus{W}{W_1})\oplus L^s_3(\brelcatrminus{W}{W_2})\]
from Lemma \ref{two}, where $W_2 <W$ and $W_1 <W$.
If $\dim W \equiv 0\Mod{4}$, we get an injection
\[L^s_3(\brelCatrminus{W\times U}) \to
 L^s_3(\brelcatrminus{W\times U}{W_1})
\oplus L^s_3(\brelcatrminus{W\times U}{W_2\times U})\]
again from Lemma \ref{two}, where $U$ is any $2${--}dimensional free
$G${--}representation.
Notice that $W_2\oplus U < W$ in the ordering above.
Consider the following commutative diagram
(with $\dim W\equiv 2\Mod{4}$):
\[
\xymatrix{
L^h_0(\ZG) \ar[rr]^{\trf{W_2\times \rminus}\ }
\ar[drr]_{\trf{W\times\rminus}}&&
L^{c_1}_1(\bcatrminus{W_2})\ar[r]\ar[d]^{\trf{W_1}}&
L^s_1(\brelCatrminus{W_2})\ar[d]^{\trf{W_1}}\\
&&L^{c_3}_3(\bcatrminus{W})\ar[r]&
L^s_3(\brelcatrminus{W}{W_1})
}
\]
This diagram, and our inductive assumption shows that
the image of  $\trf{W\times \rminus}(\sigma)$  is zero in
$L^s_3(\brelcatrminus{W}{W_1})$.
Similarly, by reversing the roles of $W_1$ and $W_2$ in the diagram,
we see that the image of $\trf{W\times \rminus}(\sigma)$  is zero in
$L^s_3(\brelcatrminus{W}{W_2})$.
Therefore the image of $\trf{W\times \rminus}(\sigma)$  is zero in
$L^s_3(\brelCatrminus{W})$.
If $\dim W \equiv 0\Mod{4}$, we replace $W$ by $W\oplus U$
and make a similar argument, using the observation that
$W_2\oplus U < W$ to justify the inductive step.
\end{proof}

\begin{proof}[The proof of Theorem \ref{one}]
Part (i) follows
from Corollary \ref{twocor}:  we know that the image of our obstruction
 $\trf{W\times\rminus}(\sigma)$ is zero ``away from $\rminus$",
so comes from a torsion element in $\hat\sigma \in L_{2k+1}^h(\bCatrminus)$
(by Lemma \ref{fourtors}, $ L_{2k+2}(\brelCatrminus{W})$ is torsion).

For part (ii) observe that
 $L^h_{2k+1}(\brelcat{\rminus}{0})=L^p_{2k}(\ZH)$ is
torsion{--}free (except for the Arf invariant $\cy 2$ which is
detected by the trivial group, so does not
matter). Hence  the image of $\trf{W\times\rminus}(\sigma)$
vanishes in $L^h_{2k+1}(\brelcat{W\times \rminus}{0})$.
The exact sequence
\[ \dots\to L^{\bu}_{2k+1}(\ZG,w) \to L^h_{2k+1}(\bcatrminus{W}) \to
L^h_{2k+1}(\brelcat{W\times\rminus}{0}),\]
with $\bU=\ker(\widetilde K_0(\ZG) \to \widetilde
K_0(\bcatrminus{W}))$, now shows that there exists a torsion element
\[\hat\sigma \in  \Image(L^{\bu}_{2k+1}(\ZG,w) \to L^p_{2k+1}(\ZG,w))\]
with $c_{\ast}(\hat\sigma) = \trf{W\times \rminus}(\sigma) \in
L^p_{2k+1}(\bcatrminus{W})$.
\end{proof}
\begin{remark}\label{sharper}
In Theorem \ref{sixteenfive} we use the sharper result
$\hat\sigma \in  L^{\bu}_{2k+1}(\ZG,w)$ to
improve our unstable similarity results
by removing extra $\rminus$ factors.
The $L$-group decoration $\bU$ is independent of $W$, since
$K_0(\bCatrminus) \to K_0(\bcatrminus{W}$ is injective
by Corollary \ref{kersimple}, and hence
$\bU=\Image(\widetilde K_0(\ZH)\to \widetilde K_0(\ZG))$.
\end{remark}

\section{ The proof of Theorem C}
 \label{C}
By replacing $V_i$ by $V_i\oplus U$ if necessary,
where $U$ is free and
$2${--}dimensional, we may assume that
$\dim V_i\oplus W\equiv 0\Mod{4}$.
This uses the $s${--}normal cobordism condition.
By the top component argument
 (see \pone, Theorem \ref{twopoweriso})
and Proposition \ref{ninethree},
we may also assume that $W$ has
only $2${--}power isotropy.
We have a commutative diagram analogous to 
\pone\ (\ref{mainda}).

\eqncount
\begin{equation} \label{maindb}
\begin{matrix}
\xymatrix@!C@C-6pc{
H^1(\widetilde K_0(\bcatrminus{W}))    \ar[dr] \ar@/^2pc/[rr]^{\delta^2}      &&
H^1(\wh(\bcatrminus{W})                \ar[dr] \ar@/^2pc/[rr]                 &&
L_0^s(\bcatrminus{W})                                                         \\
&  L_1^h(\bcatrminus{W})               \ar[dr] \ar[ur]                        &&
H^1(\Delta_{W\times\rminus})           \ar[dr] \ar[ur]                        \\
L_1^s(\bcatrminus{W})                  \ar[ur] \ar@/_2pc/[rr]_{}              &&
L_1^p(\bcatrminus{W})                  \ar[ur] \ar@/_2pc/[rr]                 &&
H^0(\widetilde K_0(\bcatrminus{W}))
}\end{matrix}
\end{equation}

\noindent where $H^1(\Delta_{W\times\rminus})$ denotes the relative group
of the double coboundary map.

This diagram for $W$ can
be compared with the one for $W_1=\Wmax$
via the inclusion maps, and we
see that the $K$-theory terms map isomorphically
by Lemma \ref{five}.
Now by Theorem \ref{one} there exists a torsion element
 \[\hat\sigma \in \Image(L^p_1(\bCatrminus)
 \to L^p_1(\bcatrminus{W_1}))\]
 which hits $\trf{W\times\rminus}$ under the cone point inclusion.
In particular, $\hat\sigma$ vanishes ``away from $\rminus$".
The main step in the proof of Theorem C is to show that the torsion
subgroup of $L^p_1(\bcatrminus{W_1})$ essentially injects into the
relative group $H^1(\Delta_{W_1\times\rminus})$ of the double
coboundary.

We have the exact sequence
\[ L^s_2(\brelCatrminus{W_1})
\to L^Y_1(\bCatrminus) \to L^s_1(\bcatrminus{W_1})
\to L^s_1(\brelCatrminus{W_1})\]
and $Y = \ker(\wh(\bCatrminus) \to \wh(\bcatrminus{W_1})$
is zero, by Corollary \ref{kersimple}.
But the first term $ L^s_2(\brelCatrminus{W_1})$ is a torsion group, 
by Lemma \ref{fourtors},
and the next term in the sequence
 $L^{\bY}_1(\bCatrminus)=L^s_1(\bCatrminus)$
 is torsion{--}free by \pone, Lemma \ref{lowA} so the middle map
$L^{\bY}_1(\bCatrminus) \rightarrowtail L^s_1(\bcatrminus{W_1})$
is an injection.

Since $L^s_1(\brelCatrminus{W_1}) =\cy 2$ by Remark \ref{fourex} we conclude
that the previous group $ L^s_1(\bcatrminus{W_1})$
is torsion{--}free modulo this $\cy 2$, which injects into
$\Lminus{1}_1(\brelCatrminus{W_1})
= \cy 2 \oplus H^1(\Kminus{1}(\ZK))$.
Since $\hat\sigma$ vanishes away from $\rminus$, this extra $\cy 2$
may be ignored.

By substituting this computation into the exact sequence
\[
L_1^s(\bcatrminus{W_1})\to L_1^p(\bcatrminus{W_1}) \to
H^1(\Delta_{W_1\times\rminus})\]
of  (\ref{maindb}), we see that the torsion subgroup of
 \[\ker (L_1^p(\bcatrminus{W_1}) \to L_1^p(\brelCatrminus{W_1})) \]
injects into the relative group
$H^1(\Delta_{W_1\times\rminus})$ of the double coboundary.
Therefore $\hat\sigma \in L_1^p(\bcatrminus{W_1})$ vanishes
if and only if
the element
 $\left\{ \Delta(V_1)/\Delta(V_2) \right\}$
is in the image of the double coboundary
\[\delta^2\colon
H^1(\widetilde K_0(\ccat{\Wmax\times \rminus,G}{\bz})) \to
H^1(\wh(\ccat{\Wmax\times \rminus,G}{\bz}))\  \]
This completes the proof of Theorem C.

\medskip
A very similar argument can be used to give an
inductive  criterion
for non{--}linear similarity without the $\rplus$ summand
(generalizing Theorem A).
In the statement we will use the analogue to
$\ck =\ker(\widetilde K_0(\ZH) \to \widetilde K_0(\ZG))$,
 namely
\[\ck_W = \ker(\wh(\brelcat{W\times\rminus}{0}) \to
\widetilde K_0(\ZG)).\]
 By Corollary \ref{kersimple}
 \[\ker(\wh(\ZG) \to 
\wh(\ccat{W\times \rminus,G}{\bz}) \cong  \Image (\wh(\ZH) \to \wh(\ZG)) \cong \wh(\ZH)\]
so we have a short exact sequence
$$0 \to \wh(\ZG)/\wh(\ZH) \to \wh(\ccat{W\times \rminus,G}{\bz})
\to \ck_W \to 0\ .$$
It follows from our $K${--}theory calculations that
$\ck_W = \ck_{\Wmax}$ and that
$\ck_W = \ck$ whenever
$\Kminus{1}(\ZK) =0$ for all $K\in \Iso(W)$.
\begin{theorem}\label{Cprime}
Let $V_1=t^{a_1} + \dots + t^{a_k} $
and $V_2 = t^{b_1} + \dots + t^{b_k}$ be
 free $G${--}representations.
Let $W$ be a  complex $G${--}representation
with no $\rplus$ summands.
Then there exists a topological similarity
 $V_1\oplus W \oplus \rminus \sim_t
V_2\oplus W \oplus \rminus$ if and only if
\begin{enumerate}
\renewcommand{\labelenumi}{(\roman{enumi})}
\item $S(V_1)$ is $s${--}normally cobordant to $S(V_2)$,
\item $\Res_H (V_1\oplus W) \oplus \rplus
 \sim_t \Res_H (V_2\oplus W)\oplus \rplus$, and
\item the element
$\left\{ \Delta(V_1)/\Delta(V_2) \right\}$
is in the image of the coboundary
\[\delta\colon
H^0(\ck_{\Wmax}) \to
H^1(\wh(\ZG^{-})/\wh(\ZH)), \]
where $0\subseteq \Wmax \subseteq W$ is a
complex subrepresentation  of real dimension $\leq 2$ with  maximal
isotropy group among the isotropy groups
of $W$ with 2-power index.
\end{enumerate}
\end{theorem}
\begin{proof}
In this case, we have a torsion element $\hat\sigma
\in L_1^h(\bcatrminus{W_1})$ which maps to our surgery
obstruction $\trf{W\times \rminus} $ under the cone point inclusion.
As in the proof of Theorem C, the torsion subgroup of
 $L_1^h(\bcatrminus{W_1})$
injects into 
$H^1(\wh(\bcatrminus{W_1}))\cong H^1(\wh(\bcatrminus{W_1}))$.
Therefore $\hat\sigma \in L_1^h(\bcatrminus{W_1})$ vanishes
if and only if
the element
 $\left\{ \Delta(V_1)/\Delta(V_2) \right\}$
is in the image of the coboundary
$$\delta\colon
H^0(\ck_{\Wmax})  \to
H^1(\wh(\ZG^{-})/\wh(\ZH))$$
as required.
\end{proof}

\section{Iterated $\rminus$ transfers}
\label{rminustwo}

We now apply Theorem \ref{optrick} to show
that iterated $\rminus$ transfers do not lead to any new
similarities.
\begin{theorem}\label{sixteenone}
Suppose that $W$ is a complex $G${--}representation
with no $\rplus$ summands.
Then $V_1\oplus W \oplus \Rminus{l}\oplus \rplus
\sim_t V_2\oplus W \oplus \Rminus{l}\oplus \rplus$
for $l\geqq 1$ implies
$V_1\oplus W \oplus \rminus\oplus \rplus
\sim_t V_2\oplus W \oplus \rminus\oplus \rplus$.
\end{theorem}
The first step in the proof is to show
injectivity of certain transfer maps.
For any homomorphism $w\colon G \to \orient$,
we will use the notation $(\ccat{W,G}{\bz},w)$ to
denote the antistructure where the involution is
$g \mapsto w(g)g^{-1}$ at the cone point.
The standard orientation (\ref{standard}) has $w=\det(\rho_W)$,
but we will need others in this section.
Let $\phi=\det(\rho_\rminus)$ for short, and notice that
$\phi\colon G \to \orient$ is  non{--}trivial.
\begin{lemma}\label{sixteentwo}
The transfer map
$\trfrminus\colon L^p_{2k+1}(\ZG,w) \to
 L^p_{2k+2}(\ccat{\rminus,G}{\bz},w\phi)$
is injective, where $w$ is the non{--}trivial orientation.
\end{lemma}
\begin{proof}
Since
$L^p_{2k+2}(\ccat{\rminus,G}{\bz},w)\cong
 L^p_{2k+2}(\rel{\ZH}{\ZG})$,
we just have to check that  the previous map in the $L^p$ version of (\ref{bdLNseq})
is zero (using the tables in \cite{ht1}).
\end{proof}
Next, a similar result for the $\Rminus{2}$ transfer.
\begin{proposition}\label{sixteenthree}
The bounded transfer
\[\trf{\Rminus{2}}\colon L^p_{2k+1}(\ZG,w) \to
 L^p_{2k+3}(\ccat{\Rminus{2},G}{\bz},w)\]
is injective, where $w$ is the non{--}trivial orientation.
\end{proposition}
\begin{proof}
We can relate the iterated $\rminus$ transfer
$\trf{1}\circ \trf{1}$ to
$\trf{2}=\trf{\Rminus{2}}$ by means of the braid diagram:
\eqncount
\begin{equation}
\begin{matrix}
\xymatrix@!C@C-5.5pc{
L_n(\ZG,w)  \ar[dr]^{\trf{1}} \ar@/^2pc/[rr]^{\trf{2}}     
&& L_{n+2}(\rel{\bz\Gamma_H}{\ZG},w)       \ar[dr]
\ar@/^2pc/[rr]   &&
LNS_{n+3}(\Phi)  \\ 
&L_{n+1}(\rel{\ZH}{\ZG},w\phi) \ar[dr] \ar[ur]^{\trf{1}}  &&
LS_{n+3}(\Phi)  \ar[dr] \ar[ur]  \\
LNS_n(\Phi)   \ar[ur] \ar@/_2pc/[rr]_{}  &&
LN_{n+3}(\rel{\ZH}{\ZG},w\phi) \ar[ur] \ar@/_2pc/[rr]                         
&& L_{n+3}(\ZG,w) }\end{matrix}
\end{equation}

\noindent where the new groups $LNS_{\ast}(\Phi)$ are the
relative groups of the transfer
\[\trf{1}=\trfrminus\colon L^p_k(\bCatrminus,w\phi) \to
 L^p_{k+1}(\bcat{\Rminus{2}},w)\ .\]
The diagram $\Phi$ of groups (as in \cite[Chap.11]{w2})
contains $\Gamma_H \cong \bz\times H$ as the
fundamental group of $S(\rminus\oplus\rminus)/G$
(see \cite[5.9]{hp2}).

The groups $LNS_{\ast}(\Phi)$ have a geometric bordism
description in terms of triples $L^{n-2}\subset N^{n-1} \subset M^n$
where each manifold has fundamental group $G$, $N$ is
characteristic in $M$, and $L$ is characteristic in $N$, with
respect to the index two inclusion $H < G$.
There is also an algebraic description for $LNS_{\ast}(\Phi)$
in terms of ``twisted antistructures" \cite{htw3} as for the other
groups $LN(\rel{\ZH}{\ZG},w\phi)$ \cite[12C]{w2} and
for $LS(\Phi,w)$ \cite[7.8.12]{ra8}.
Substituting these descriptions into the braid gives:

\eqncount
\begin{equation}
\begin{matrix}
\xymatrix@!C@C-6pc{
L_{2k+1}(\ZG,\alpha,u)   \ar[dr]^{\trf{1}} \ar@/^2pc/[rr]^{\trf{2}}      &&
L_{{2k+1}}(\rel{\bz\Gamma_H}{\ZG},\tilde\alpha,\tilde u)      \ar[dr] \ar@/^2pc/[rr]     &&
L_{{2k}}(\rel{\ZH}{\bz\Gamma_H},\tilde \alpha,\tilde u)       \\
&   L_{2k+1}(\rel{\ZH}{\ZG},\alpha,u)                        \ar[dr] \ar[ur]^{\trf{1}}    &&
L_{{2k}}(\bz\Gamma_H,\tilde\alpha,\tilde u)               \ar[dr] \ar[ur]     \\
L_{2k+1}(\rel{\ZH}{\bz\Gamma_H},\tilde \alpha,\tilde u)       \ar[ur] \ar@/_2pc/[rr]_{}    &&
L_{2k}(\ZH,\alpha,u)                                   \ar[ur] \ar@/_2pc/[rr]     &&
L_{2k}(\ZG,\alpha,u)
}\end{matrix}
\end{equation}

\noindent The antistructure $(\ZG,\alpha,u)$ is the twisted antistructure
obtained by scaling with an element  $a\in G - H$.
The antistructure $(\bz\Gamma_H,\tilde\alpha,\tilde u)$ is
the one defined by Ranicki \cite[p.805]{ra8}, then
scaled by $\tilde a \in \Gamma_H$, where $\tilde a$
maps to $a$ under the projection $\Gamma_H \to G$.
Since $\Gamma_H \cong \bz\times H$, we have an exact
sequence (see \cite[Theorem 4.1]{milgran2})
\[\to L^p_n(\ZH,\alpha,u) \to
 L^p_n(\bz\Gamma_H,\tilde\alpha,\tilde u)\to
 \Lminus{1}_{n-1}(\ZH,\alpha,u)\RA{1-w(a)}
L^p_{n-1}(\ZH,\alpha,u)\]
and it follows that
$L^p_n(\rel{\ZH}{\bz\Gamma_H},\tilde \alpha,\tilde u)
\cong \Lminus{1}_{n-1}(\ZH,\alpha,u)$.
It is not difficult to see that
$\Lminus{1}_{2k}(\ZH,\alpha,u)$ is torsion{--}free (except for
the Arf invariant summand)
by a similar argument to \pone, Theorem \ref{Lfacts},
using the $L^p$ to $\Lminus{1}$ Ranicki-Rothenberg sequence.
We  first check that
 $L^p_{2k}(\ZH,\alpha,u)$ is torsion{--}free (again except
for the Arf invariant summand)  from the tables
in \cite[14.21]{ht1}.

Now observe that $L^p_n(\ZG,\alpha,u)\cong L^p_n(\ZG,w)$,
and
\[L^p_{2k+1}(\rel{\ZH}{\ZG},\alpha,u) \cong L^p_{2k+2}(\rel{\ZH}{\ZG})\ .\]
The transfer
 $\trfrminus\colon L^p_{2k+1}(\ZG,w) \to  L^p_{2k+2}(\rel{\ZH}{\ZG})$
is  injective by Lemma \ref{sixteentwo}, and
since $w(a)=-1$, the map
\[ \Lminus{1}_{2k}(\ZH,\alpha,u)\RA{1-w(a)}
L^p_{2k}(\ZH,\alpha,u)\] is also injective (except for the Arf invariant summand).
Therefore $ L^p_{2k+1}(\ZG,w)$ must
inject into the relative group
$L^p_{2k+1}(\rel{\bz\Gamma_H}{\ZG},\tilde\alpha,\tilde u)\cong
L^p_{2k+3}(\bcat{\Rminus{2}},w)$.
\end{proof}
The final step in the proof is to consider the following
 commutative  square of spectra:
\[
\xymatrix{\BL (\ccat{pt}{\ZG},w)
\ar[r]^{c_{*}\ }\ar[d]_{\trfrminus}
&\BL (\bcat{V},w)
\ar[d]^{\trfrminus}\\
\ \BL (\ccat{\rminus,G}{\bz},w\phi)\ar[r]^{c_{*}\ }
&\ \BL (\bcat{V\times \rminus},w\phi)
}
\]
for any $G$-representation $V$ with no $\rplus$ summands.
\begin{lemma}\label{sixteenfour}
This is a pull{--}back square of $\Lminus{\infty}$  spectra.
\end{lemma}
\begin{proof}
The fibres of the vertical $\rminus$ transfers are
\[\mathbb{LN}(\bCatrminus, w\phi) \cong
\BL(\ccat{pt}{\ZH ,\alpha,u})\ ,\] and
$\mathbb{LN}(\bcatrminus{V},w\phi) \cong
\BL(\ccat{V,H}{\bz},\alpha,u)$ respectively by
(\ref{bdLNseq}),
and the fibre of the cone point inclusion
\[\BL(\ccat{pt}{\ZH},\alpha,u) \to \BL(\ccat{V,H}{\bz},\alpha,u)\]
is  $\BL(\Ccat{V,H}{0}{\bz},\alpha,u)$.
But Theorem \ref{optrick} shows that
$\BL(\Ccat{V,H}{0}{\bz},\alpha,u)$ is contractible, and
therefore the cone point inclusion is a homotopy equivalence.
It follows that the fibres of
the horizontal maps are also homotopy equivalent.
\end{proof}
\begin{proof}[The proof of Theorem \ref{sixteenone}]
It is enough to prove that $\trf{W\times\Rminus{l+1}}(\sigma)=0$, for
$l\geqq 1$, implies $\trf{W\times\rminus}(\sigma)=0$ in the top
component of $L^p_{2k+1}(\bcatrminus{W})$. We may assume that
$\Iso(W)$ only contains subgroups of $2${--}power index, and let $\dim
W =2k$. We may also assume that $l$ is even, by crossing
with another $\rminus$ if necessary.
Now by Theorem C, $\trf{W\times\Rminus{2s+1}}(\sigma)=0$
implies $\trf{W\times\Rminus{2s-1}}(\sigma)=0$,
provided that $s\geqq 2$.
It therefore remains to study  $l=2$.

The pullback squares provided by Lemma \ref{sixteenfour} can
be combined as follows. Consider the diagram of 
$\Lminus{\infty}$-spectra

$$\xymatrix{
\BL  (\ccat{pt}{\ZG},w)\ar[r]^{c_{*}\ }\ar[d]_{\trfrminus}&\BL  (\ccat{\rminus,G}{\bz},w\phi)
\ar[r]^{c_{*}\ }\ar[d]_{\trfrminus}&
\BL  (\ccat{W\times\Rminus{2},G}{\bz},w)\ar[d]_{\trfrminus}\cr
\BL  (\ccat{\rminus,G}{\bz},w\phi)\ar[r]^{c_{*}\ }&\BL  (\ccat{\Rminus{2},G}{\bz},w)\ar[r]^{c_{*}\ }&
\BL  (\ccat{W\times\Rminus{3},G}{\bz},w\phi)
}$$
whose outer square ($V=W\times\Rminus{2}$)
 and the left-hand square ($V=\rminus$) are both pullback squares, and hence so is the right-hand
square. Next consider the diagram of $\Lminus{\infty}$-spectra
$$\xymatrix{
\BL  (\ccat{pt}{\ZG},w)\ar[r]^{c_{*}\ }\ar[d]_{\trfrminus}&
\ \BL  (\ccat{W\times\rminus,G}{\bz},w\phi)\ar[d]_{\trfrminus}\cr
\BL  (\ccat{\rminus,G}{\bz},w\phi)
\ar[r]^{c_{*}\ }\ar[d]_{\trfrminus}&
\ \BL  (\ccat{W\times\Rminus{2},G}{\bz},w)\ar[d]_{\trfrminus}\cr
\BL  (\ccat{\Rminus{2},G}{\bz},w)\ar[r]^{c_{*}\ }&
\ \BL  (\ccat{W\times\Rminus{3},G}{\bz},w\phi)
}$$
The lower square was just shown to be a pull-back, and the upper
square is another special case of Lemma \ref{sixteenfour}
(with $V= W\times\rminus$). Therefore the outer square is
a pull-back, and this is the one used for the case $l=2$.

Now we apply homotopy groups to these pull-back squares
(using the fact that $\Lminus{\infty} \equiv \Lminus{1}$  to obtain  the lower
squares in the commutative diagram:
\[
\xymatrix{L^p_{2k+2}(\brelcat{W\times\rminus}{0})\ar[r]\ar[d]^{\cong}
&L^p_{2k+1}(\ZG,w)\ar[r]^{c_{\ast}}\ar[d]
&L^p_{2k+1} (\bcat{W\times\rminus},w)\ar[d]\\
\Lminus{1}_{2k+2}(\brelcat{W\times\rminus}{0})
\ar[r]\ar[d]^\cong
&\Lminus{1}_{2k+1}(\ZG,w)\ar[r]^{c_{\ast}}\ar[d]_{\trf{\Rminus{2}}}
&\Lminus{1}_{2k+1} (\bcat{W\times\rminus},w)\ar[d]^{\trf{\Rminus{2}}} \\
\Lminus{1}_{2k+4}(\brelcat{W\times\Rminus{3}}{0})\ar[r]
&\Lminus{1}_{2k+3}(\ccat{\Rminus{2},G}{\bz},w)\ar[r]^{c_{\ast}}&
\Lminus{1}_{2k+3}(\bcat{W\times \Rminus{3}},w)
}
\]
 The vertical maps in the top row squares are induced by the change of $K$-theory decoration. 
We need some information 
about the maps in this diagram.

\begin{enumerate}
\renewcommand{\labelenumi}{(\roman{enumi})}
\item The upper left-hand vertical
map is an isomorphism, since 
$\Kminus{1}(\brelcat{W\times\rminus}{0})=0$.

\item The lower left-hand vertical map is an isomorphism, 
by Lemma \ref{sixteenfour}.

\item The transfer map
$$\trf{\Rminus{2}}\colon L^p_{2k+1}(\ZG,w) \to
 L^p_{2k+3}(\bcat{\Rminus{2}},w)$$ is injective, from
Proposition \ref{sixteenthree}.

\item 
The map 
$$L^p_{2k+3}(\ccat{\Rminus{2},G}{\bz},w)
\to \Lminus{1}_{2k+3}(\ccat{\Rminus{2},G}{\bz},w)$$ is
injective, by \pone, Corollary \ref{rstable}.

\item The composite of the middle two vertical maps
in the diagram is injective,  by combining
parts (iii)-(iv).

\end{enumerate}

\noindent
Suppose 
 that  $\trf{W\times\Rminus{3}}(\sigma)=0$.
We have proved that the transfer map
\[\trf{\Rminus{2}}\colon L^p_{2k+1}(\ZG,w) \to
\Lminus{1}_{2k+3}(\ccat{\Rminus{2},G}{\bz},w)\]
is injective in part (v).
Since $\trf{W\times\rminus}(\sigma) = c_{\ast}(\hat\sigma)$
for some $\hat\sigma \in L^p_{2k+1}(\ZG,w)$, 
 it follows  by a diagram chase that 
 $\trf{W\times\rminus}(\sigma)=0$ and we are done.
\end{proof}
We also have a version without $\rplus$ summands.
\begin{theorem}\label{sixteenfive}
Suppose that $W$ is a complex $G${--}representation
with $W^G=0$.
Then $V_1\oplus W \oplus \Rminus{l}
\sim_t V_2\oplus W \oplus \Rminus{l}$
for $l\geqq 3$ implies
$V_1\oplus W \oplus \Rminus{3}
\sim_t V_2\oplus W \oplus \Rminus{3}$.
\end{theorem}
\begin{proof}
This follows from a similar argument, using Theorem \ref{Cprime}
instead of Theorem C.
The injectivity results of Lemma \ref{sixteentwo} and
Proposition \ref{sixteenthree} also hold for $\trf{1}$ and $\trf{2}$
on $L^{\bu}_{2k+1}(\ZG,w)$,
where $\bU=\Image(\widetilde K_0(\ZH) \to \widetilde K_0(\ZG))$
as in Remark \ref{sharper}. This $L$-group decoration  fits
in exact sequences with $L^p_{n}(\ZH)$ so the previous calculations
for injectivity apply again.
The details are left to the reader.
\end{proof}

\section{The proof of Corollary \ref{sixdim}}
\label{sixdiml}
Suppose that $G=\Cy{4q}$, $q$ odd, and that $V_1\oplus \rminus\oplus
\rplus\sim_t V_2\oplus \rminus\oplus
\rplus$ with $\dim V_i =4$.
We will use Theorem B and the assumption about odd class numbers
to prove that $V_1\cong V_2$.
This is the only case we need to discuss to prove Corollary
\ref{sixdim}. If we started instead with
$V_1\oplus \Rminus{2}\sim_t V_2\oplus \Rminus{2}$,
then stabilizing with $\rplus$ and applying Theorem \ref{sixteenone} would reduce
to the case above.

We may assume that $q > 1$ is a prime, and that the free
representations have the form
 $V_1 = t + t^i$ and $V_2=t^{1+2q} + t^{i+2q}$
where $(i,2q)=1$ (not all weights are $\equiv 1 \Mod 4$).
The Reidemeister torsion invariant is
$\Delta(V_1)/\Delta(V_2) = U_{1,i} \in \wh(\ZG)$
in the notation for units of $\ZG$ used in
\cite[p. 732]{cssw1}.
With respect to
 the involution $\Bar{\hphantom{a}}\colon t \mapsto -t^{-1}$, this element
defines a class $u_{1,i} = \{U_{1,i}\} \in H^1(\wh(\ZG^{-}))$.
Since the map $L^p_1(\ZG^{-}) \to L^p_1(\bCatrminus)$ is injective,
it is enough to show that the image of  $u_{1,i}$ is non{--}trivial
in $\Coh{1}{\Delta}$ (see diagram \pone\ (\ref{lpbraid})).

Let $A = \bz[\zeta_{4q}]$ be the ring of cyclotomic integers.,
and $B=\bz[\zeta_{2q}]$.
We will study the top component $\wh(\ZG)(q)$ by
comparing it to $K_1(\cm)$, where $\cm = A \times B\times B$
is the top component of an involution invariant maximal order
in $\QG$ containing $\ZG$.
Notice that the two copies of $B$ are
interchanged under the involution, so
$\Coh{i}{K_1(\cm)} = \Coh{i}{\Units{A}}$.
We have the exact sequence
\[0 \to \wh(\ZG)(q) \to K_1(\cm)
\to K_1( \cmlocal)/\wh(\zlocal G)(q) \to D(\ZG)(q)\to 0\ .\]
where $D(\ZG) =\ker(\widetilde K_0(\ZG) \to \widetilde K_0({\cm}))$.
Note that $H^i(D(\ZG)) =H^i(\widetilde K_0(\ZG))$
since $A$ has odd class number.

Let $i_{\ast}\colon \Coh{1}{\wh(\ZG^{-})(q)} \to \Coh{1}{\Units{A}}$
denote the map on the Tate cohomology
induced by the inclusion $ i\colon\wh(\ZG)(q)\to \Units{A}$.
\begin{lemma}\label{seventeenone}
For  $(j,2q)=1$ and $j\equiv 1\Mod{4}$,
the image $i_{\ast} (u_{1,j}) = \langle -1\rangle \in H^1(\Units{A})$.
\end{lemma}
\begin{proof}
Let $\gamma_j = \frac{t^j-1}{t-1} \in \Units{A}$, and compute
$\bar \gamma_j/\gamma_j$ with respect to the
non{--}oriented involution.
We get $\bar \gamma_j/\gamma_j = u_{1,j}/u_{1,1}$,
so $i_{\ast}(u_{1,j}) = i_{\ast}(u_{1,1})$ for all $i$.
Now let $v = \frac{t+1}{t-1}$.
Since $t=-t^{1+2q}$, this is a cyclotomic
unit in $A$ with $v\bar v = -1$.
But $v/\bar v = -u_{1,1}$, so $i_{\ast} (u_{1,j}) = \langle -1\rangle$
in $H^1(\Units{A})$.
\end{proof}
Next we need a computation:
\begin{lemma}
 $ 0 \neq \langle -1\rangle \in H^1(\Units{A})$.
\end{lemma}
\begin{proof}
Let $E = \bq(\zeta_{4q})$ be the quotient field of $A$.
Since $A$ has odd class number,
$\Coh{i}{\Units{E}/\Units{A}} \cong
\Coh{i}{\localunits{E}{}{}/\localunits{A}{}{}}$.
Next, observe that the extension $E/F$ is evenly ramified at $q$,
where $F$ is the fixed field of the involution
$\Bar{\hphantom{a}}\colon E \to E$.
It follows that the map $\Coh{0}{\localunits{E}{}{q}}
\to \Coh{0}{\localunits{E}{}{q}/\localunits{A}{}{q}}$
is zero.
We finish by considering the commutative diagram
\[
\xymatrix{\Coh{0}{\Units{E}}\ar[r]\ar[d]&
\Coh{0}{\Units{E}/\Units{A}}\ar[r]^{\delta}\ar[d]^{\approx}&
\Coh{1}{\Units{A}}\ar[d]\\
\Coh{0}{\localunits{E}{}{}}\ar[r]&
\Coh{0}{\localunits{E}{}{}/\localunits{A}{}{}}\ \ar@{>->}[r]&
\Coh{1}{\localunits{A}{}{}}
}
\]
and the prime element $w=\zeta_q - \zeta_q^{-1} \in E$
over $q$.
Since $\bar w/w = -1$,  $\omega$  defines an element in
$\Coh{0}{\Units{E}/\Units{A}}$ and
$\delta(w) = \langle -1\rangle \in H^1(\Units{A})$.
On the other hand,  $w$ maps under the middle isomorphism
to a non{--}zero element in
$\Coh{0}{\localunits{E}{}{q}/\localunits{A}{}{q}}$ so
$\delta(w) \neq 0$.
\end{proof}
Now we can complete the proof of Corollary \ref{sixdim} by considering
the commutative braid:

\eqncount
\begin{equation}
\begin{matrix}
\xymatrix@!C@C-4pc{
\Coh{0}{\Units{A}}                                                         
\ar[dr]           \ar@/^2pc/[rr]                &&
H^0\left (\frac{K_1(\cmlocal)}{\wh(\zlocal G)(q)}\right )       \ar[dr]
\ar@/^2pc/[rr]                          &&
\Coh{0}{\widetilde K_0(\ZG^{-})(q)}                                            
                                                \\
&  H^0\left (\frac{K_1(\cm)}{\wh(\ZG)(q)}\right )                  \ar[dr]
\ar[ur]                                 &&
\Coh{1}{\Delta}                                                            
\ar[dr] \ar[ur]                    \\
\Coh{1}{\widetilde K_0(\ZG^{-})(q)}                                             \ar[ur]
\ar@/_2pc/[rr]^{\delta^2}               &&
\Coh{1}{\wh(\ZG^{-})(q)}                                             
\ar[ur] \ar@/_2pc/[rr]^{i_*}                    &&
\Coh{1}{\Units{A}}
}\end{matrix}
\end{equation}

\noindent containing the double coboundary and its relative group
$\Coh{1}{\Delta}$.
By commutativity of the lower right triangle, the image of
our obstruction
element $u_{1,j}$ is non{--}zero in $\Coh{1}{\Units{A}}$
and therefore non{--}zero in $\Coh{1}{\Delta}$.

\section{The normal invariant}
\label{normalinvariant}

In this section we collect some results about
 normal cobordisms of lens spaces with cyclic fundamental group
$G$. Recall that a necessary condition for the existence of
a similarity $V_1\oplus W \sim_t V_2\oplus W$ is that $S(V_1)$ and
$S(V_2)$ must be $s$-normally cobordant.
In \cite[\S 1]{cssww1} it is asserted
that their formula $(A')$ gives necessary and sufficient
conditions for homotopy equivalent lens spaces
to be $s${--}normally cobordant, when $G = \Cy{2^r}$.
However in \cite[1.3]{yo1} it is shown that the given conditions $(A')$
are sufficient but not necessary.
We only use the sufficiency here, and 
 study the subgroups
$\rtildefree{G}$ and $\Rtildefree_h(G)$ of $R(G)$
defined in Section \ref{splitting}.
Recall that 
 $\rtildefree{G} = \ker(\Res\colon\rfree{G} \to \rfree{\Godd})$.
\begin{lemma}\label{twentyone}
Let $G = \Cy{ 2^r q}$ be a finite cyclic group, and
$G_2$ be the 2{--}Sylow subgroup.
If $q >1$ there is an exact sequence
\[0 \to \cy 2 \to  \rtildefree{G}/ \Rtildefree_h(G) \to
\Rtildefree(G_2)/ \Rtildefree_h(G_2) \to 0\]
given by the restriction map $\Res _{G_2}$.
The kernel is generated by any element
$\alpha \in  \rtildefree{G}$ with  $k${--}invariant
$k(\alpha) \in \Units{(\cy{|G|})}/\{\pm 1\}$
in the coset   $ k \equiv 1\Mod{2^r}$ and
$k \equiv -1\Mod{q}$.
\end{lemma}
\begin{proof}
There is a short exact sequence
\[0 \to  \Rtildefree_h(G) \to  \rtildefree{G} \to
 \Units{(\cy{ 2^rq})} /\{ \pm 1\}\]
given by the $k${--}invariant.
Moreover,
the  $k${--}invariants of elements in  $\rtildefree{G} $
lie in
\[\ker \bigl (  \Units{(\cy{2^rq})}/\{ \pm 1\}
 \to  \Units{(\cy{q})} /\{ \pm 1\}\bigr )
\cong \{lq+1\Vertical l \in \bz\}
\subseteq   \Units{(\cy{2^rq})}\]
which injects  under reduction $\Mod{2^r}$ into
$\Units{(\cy{2^r})}$ if $q >1$.
Restriction of the $k${--}invariant to $G_2$ detects only
its value in  $\Units{(\cy{2^r})} /\{ \pm 1\}$, so the
kernel has order 2.
Since the restriction map
$\Rtildefree{G} \to  \Rtildefree(G_2)$
is surjective, we have the required exact sequence.
\end{proof}
The situation for the normal invariant is simpler.
Recall that we write $V \bumpeq V'$  or $(V-V')\bumpeq 0$
if  there exists a homotopy equivalence
$f\colon S(V)/G \to S(V')/G$ of lens
spaces, such that $f$ is $s${--}normally cobordant to the identity.
\begin{lemma}\label{twentytwo}
Let $G = \Cy{2^rq}$ be a finite cyclic group.
Then
\[
\ker \bigl (\Res\colon \Rfree(G) \to
 \Rfree(\Cy{ q})\oplus \Rfree(\Cy{2^r})\bigr )
\subseteq \Rfree_n(G)\ .
\]
\end{lemma}
\begin{proof}
We may assume that $r\geq 2$, since
$\Res\colon  \Rfree(\Cy{2q}) \to  \Rfree(\Cy{q})$
is an isomorphism, and consider an element
\[\sum t^{a_i} - \sum t^{b_i} \in \ker \bigl (\rfree{G} \to
\Rfree(\Cy{q})\oplus\Rfree(\Cy{2^r})\bigr ).\]
With a suitable ordering of the indices, $a_i \equiv b_{i}
\Mod{2^r}$ and
$a_i \equiv b_{\tau(i)}\Mod{q}$ for some
permutation $\tau$.
Since the normal invariants of lens spaces are detected
by a cohomology theory, we can check the condition
at each Sylow subgroup separately.
It follows that $\sum t^{a_i} \bumpeq \sum t^{b_i}$.
\end{proof}
\begin{example}
For $G=\Cy{24}$ we have $\Rtildefree(G) = \{t-t^5, t-t^7, t-t^{11}\}$
and the subgroup 
$\ker\Res_{\Cy{8}}\cap \Rtildefree(G) = \{ t-t^7, t^5 - t^{11}\}$.

\end{example}

\begin{lemma}\label{twentythree}
Let $G = \Cy{2^rq}$ be a finite cyclic group, and
$G_2$ be the 2{--}Sylow subgroup.
If $r \geq 2$ there is an isomorphism
\[ \Rtildefree_h(G)/ \Rtildefree_n(G) \to
\Rtildefree_h(G_2)/ \Rtildefree_n(G_2) \]
given by the restriction map $\Res _{G_2}$.
\end{lemma}
\begin{proof}
We first remark that if
$W = t^{a_1} + \dots + t^{a_n}$ is a free
 $G_2${--}representation (here we use the assumption
that $r\geq 2$), then by choosing
integers $b_i \equiv a_i\Mod{2^r}$ and
$b_i \equiv 1\Mod{q}$ we obtain a free
$G${--}representation
$V = t^{b_1} + \dots + t^{b_n}$  with
$\Res _{G_2}(V) = W$.
It follows that
\[ \Res _{G_2}\colon \Rtildefree_h(G) \to
\Rtildefree_h(G_2) \]
is surjective, and therefore
\[ \Res _{G_2}\colon \Rtildefree_h(G)/ \Rtildefree_n(G)) \to
\Rtildefree_h(G_2)/ \Rtildefree_n(G_2) \]
is also surjective.

Now suppose $\alpha \in \Rtildefree_h(G)$
and $\Res _{G_2}(\alpha) \in \Rtildefree_n(G_2)$.
This means that $\Res _{G_2}(\alpha) = (W -W')$  for some
free ${G_2}${--}representation $W$, $W'$ such
that $W \bumpeq W'$.
By the construction of the last paragraph, we can
find free $G${--}representations $V$, $V'$
such that
(i) $\Res _{G_2}(V) = W$ and $\Res _{G_2}(V') = W'$, and
(ii) $\alpha'=(V - V') \in \Rtildefree_h(G)$.
It follows that
$V\bumpeq V'$ and so $\alpha'  \in
\Rtildefree_n(G))$, and
\[\alpha - \alpha' \in
\ker \bigl (\Rfree_h(G) \to
 \Rfree_h(\Cy{q})\oplus \Rfree_h(\Cy{2^r})\bigr ).\]
By Lemma \ref{twentytwo}, $\alpha - \alpha' \in \Rtildefree_n(G)$
and so $\alpha \in \Rtildefree_n(G)$.
\end{proof}
Our final result is a step towards determining $R_h(G)/R_n(G)$
more explicitly.
\begin{lemma}\label{twentyfour}
Let $ b_{i,s}=(t^i - t^{2^{r-s}q-i}) \in \rtildefree{G}$
for $1\leq i <2^{r-s-1}q$, $1\leq s \leq r-1$, and $ (i,2q) =1$.
Let $l(i,s) = s$ for $1\leq s \leq r-2$, and
$l(i,r-1)$ the order of $i/(2q-i) \in \Units{(\cy{2^r})}$.
Then $b_{i,s}$ is an element of order $2^{l(i,s)}$
in  $R(G)/R_h(G)$, and
 $2^{l(i,s)}\cdot b_{i,s}\bumpeq 0$.
\end{lemma}
\begin{proof}
If $1\leq s\leq r-2$ then $k(b_{i,s}) = i/(2^{r-s}q-i) \in
 \Units{(\cy{2^rq})}/\{ \pm 1\}$ is congruent to $-1\Mod{4q}$.
By Lemma \ref{twentyone}  the linear span of these elements
in $R(G)/R_h(G)$ injects into $R({G_2})/R_h({G_2})$, where
${G_2}$ is the 2{--}Sylow subgroup.
By Lemma \ref{twentythree} the normal invariant  is also detected
by restriction to ${G_2}$, so it is enough to prove the
assertions about these elements ($s\leq r-2$) when $G=\Cy{2^r}$.

For the first part we must show that the $k${--}invariant
of $b_{i,s}$ has order $2^{s}$.
We will use the expression
\[\nu_2(r!) = r-\alpha_2(r),\]
where $\alpha_2(r)$ is the number of non-zero  coefficients
in the 2-adic expansion of $r$, for the 2-adic valuation of $r$.
Then
\[\nu_2(\binom{2^s}{k}) = s-\nu_2(k)\]
for $1\leq k\leq 2^s$.
These formulas and the binomial expansion show that
$k(b_{i,s})^{2^{s}} \equiv 1 \Mod{2^r}$ for $1\leq s\leq r-2$.

Next we consider the normal invariant.
For $1\leq s\leq r-2$ we will take $q=1$ and
 apply the criterion $(A')$ of \cite{cssww1}
to show that $2^{s}\cdot b_{i,s}\approx 0$.
This amounts to a re-labelling of our original
elements $b_{i,s}$ without changing the order
of their $k${--}invariants.

We must now compute the elementary
symmetric functions $\sigma_k(2^s \cdot i^2)$
and $\sigma_k(2^s \cdot (2^{r-s}-i)^2)$,
where the notation $2^s \cdot i^2$ means that the weight
$i^2$ is repeated $2^s$ times in the symmetric function.
The formula $(A')$ is
\[\sigma_k(2^s \cdot (2^{r-s}+i)^2) -
\sigma_k(2^s \cdot i^2) \equiv 2\bigl ( (2^{r-s}-i)^{2^s}
- i^{2^s}\bigr )\binom{2^s-1}{k-1} \Mod{2^{r+3}}\]
so on the right-hand side we have
\[
\begin{cases} +2^{r+1} \Mod{2^{r+3}}& \text{if $k$ odd}\\
 -2^{r+1}\Mod{2^{r+3}} &\text{if $k$ even}
\end{cases}
\]
The formula $(A')$ assumes that the weights are congruent
to $1 \Mod{4}$, or in our case $i \equiv 1\Mod{4}$.
Therefore, if $i \equiv 3\Mod{4}$, we must use the equivalent
weights $-i$ and $-(2^{r-s}-i)$.

To compute the left-hand side we use the Newton polynomials
$s_k$ and their expressions in term of elementary symmetric
functions. We need the property $s_k(2^s\cdot i^2) = 2^s  i^{2k}$
and the coefficient of $\sigma_k$ in $s_k$
which is $(-1)^{k+1}k$.
By induction, we see that the left-hand side is just
\[\frac{(-1)^{k+1}}{k} (s_k(2^s \cdot (2^{r-s}-i)^2) -
s_k(2^s \cdot i^2)\bigr ) \equiv \frac{(-1)^{k+1}2^s}{k}
\bigl ( (2^{r-s}-i)^{2k} - i^{2k}   )\Mod{2^{r+3}}\]
But by writing
\[(2^{r-s}-i)^{2k} - i^{2k}  = (2^{r-s}\theta + i^2)^k - i^{2k}\]
where $\theta = 2^{r-s} - 2i \equiv 2 \Mod{4}$, our expression
becomes
\[\frac{(-1)^{k+1}2^s}{k}(2^{r-s}k\theta)
\equiv (-1)^{k+1}2^{r+1}\Mod{2^{r+3}}\]

If $s = r-1$, then $k(b_{i,s}) = i/(2q-i) \equiv -1 \Mod{q}$,
and $k(b_{i,s}) \equiv 1 \Mod{2^r}$ whenever $i \equiv q
\Mod{2^{r-1}}$. This is for example always the case
for $G = \Cy{4q}$.
Such  elements
$b_{i,r-1}$ lie in the kernel of the restriction map to $R(\Cy{2^r})$.
Moreover, if $b_{i,r-1}
 \in \ker\bigl (R(G) \to R(\Cy{2^r})\bigr )$, then
$2b_{i,r-1} \in R_h(G)$
and $2b_{i,r-1}\bumpeq 0$ by Lemma \ref{twentytwo}.
Otherwise, the order of $b_{i,r-1}\in R(G)/R_h(G)$ is
the same as the order of  its restriction
 $\Res (b_{i,r-1}) $ to the 2{--}Sylow subgroup,
and $\Res (b_{i,r-1}) = \pm \Res (b_{j,s})$ for
$l(i,r-1)=s \leq r-2$ and some $j$.
\end{proof}
\begin{remark} As pointed out by the referee, Lemma 12.1 and
Lemma 12.4 together give a short exact sequence
$$0 \to \cy 2 \to \Rtildefree(G)/ \Rtildefree_n(G) 
\to  \Rtildefree(G_2)/ \Rtildefree_n(G_2) \to 0, $$ assuming
that $G$ is not a cyclic $2$-group. It is not difficult to see from
the proof of Theorem E (given in \pone), that for $G_2=\Cy{2^r}$, 
and $r\geqq 4$, the term
$ \Rtildefree(G_2)/ \Rtildefree_n(G_2)$ is the quotient
of $\rtildetopfree{\Cy{2^r}}$ by the subgroup $\langle \alpha_1 +\beta_1\rangle$.
\end{remark}

\section{The proof of Theorem D}
\label{sD}
We first summarize our information about $\rtop{G}$, obtained
by putting together results from previous sections.
If $\alpha \in \Rtildefree_{h,\Top}(G)$, we define 
the \emph{normal invariant order} of $\alpha$ to be the minimal
$2${--}power such that $2^t\alpha \in \Rtildefree_{n,\Top}(G)$.
\begin{theorem}\label{oneD}
Let $G=\Cy{2^rq}$, with $q$ odd, and $r \geqq 2$.
\begin{enumerate}
\renewcommand{\labelenumi}{(\roman{enumi})}
\item The torsion subgroup of $\rtopfree{G}$ is $\rtildetopfree{G}$.
\item The rank of $\rtopfree{G}$ is $\varphi(q)/2$ for $q >1$
 (resp. rank 1 for $q=1$), and the torsion is at most $2${--}primary.
\item The subgroup $\Rtildefree_{n,\Top}(G)$  has exponent two,
and the Galois action induced by group automorphisms is the identity.
\item For any $\alpha \in \Rtildefree_{h,\Top}(G)$,
if the normal invariant order of $\Res_H(\alpha)$ is $2^t$,
then the normal invariant order of $\alpha$ is $2^{t+1}$.
\end{enumerate}
\end{theorem}
\begin{remark}
In part (ii), $\varphi(q)$ is the Euler function.
The precise number   of $\cy 2$ summands in
$\Rtildefree_{n,\Top}(G)$ is determined by working out
the conditions in Theorem C on the basis elements of
$\Rtildefree_n(G)$.
In cases where the conditions in Theorem C can
actually be evaluated,
the structure of $\rtop{G}$ will thus be determined completely.
\end{remark}
\begin{proof}
Parts (i) and (ii) of Theorem \ref{oneD} have already been proved in
Corollary \ref{cyclicranks}, so it remains to discuss parts (iii) and (iv).
In fact, the assertion that $\Rtildefree_{n,\Top}(G)$
has exponent $2$ is an immediate consequence of
Theorem C.
To see this, suppose that $(V_1 - V_2)$ is any element
 in $\Rtildefree_{n}(G)$.
The obstruction
to the existence of a stable non{--}linear similarity
$V_1\approx_t V_2$ is determined by the class
$\{\Delta(V_1)/\Delta(V_2)\}$ in  the Tate cohomology group
$H^1(\wh(\ZG^{-})/\wh(\ZH))$, which has exponent $2$.
Since the Reidemeister torsion is multiplicative,
$\Delta (V_1 \oplus V_1) = \Delta(V_1)^2$, and
we conclude that $V_1\oplus V_1\approx_t V_2\oplus V_2$
by Theorem C.
Finally, suppose that $\alpha \in \Rtildefree_{n, \Top}(G)$ and
that $\sigma$ is a group automorphism of $G$. By induction, we
can assume that $\Res_H(\alpha - \sigma(\alpha)) = 0$.
We now apply Theorem C to $\beta=\alpha - \sigma(\alpha)$,
with  $W$  a complex
$G${--}representation such that $W^G=0$, containing all the
non-trivial irreducible representations of $G$ with isotropy of 
$2$-power index. Then $\beta$ is
detected by  the image of its Reidemeister torsion invariant
in $H^1(\Delta_{W\times\rminus})$. But by   \pone,  Lemma
\ref{galoisinvariance}, the Galois action on this group is trivial.
Therefore $\beta=\alpha - \sigma(\alpha)=0$.

For part (iv) we recall that the normal invariant for $G$ lies in
a direct sum of groups $H^{4i}(G;\bz_{(2)}) \cong \cy{2^r}$.
Since the map $\Res_H$ induces the natural projection
$\cy{2^r} \to \cy{2^{r-1}}$ on group cohomology, the result follows
from \cite[2.6]{cs1}.
\end{proof}

It follows from our results that
the structure of
$\Rtildefree_{n,\Top}(G)$ is determined
by working out the criteria of Theorem C on a basis of
$\Rtildefree_n(G)$.
Suppose that  $V_1$ is stably topologically equivalent to $V_2$.
Then there exists a similarity
\[V_1\oplus W \oplus \Rminus{l} \oplus \Rplus{s}
\sim_t V_2\oplus W \oplus \Rminus{l} \oplus \Rplus{s}\]
where $W$ has no $\rplus$ or $\rminus$ summands,
and $l,s\geqq1$.
But by \pone, Corollary \ref{rstable} we may assume that $s=1$,
and by Theorem \ref{sixteenone} that $l=1$, so we are reduced
to the situation handled by
Theorem C. The algebraic indeterminacy given there is computable,
but not very easily if the associated cyclotomic fields have
complicated ideal class groups. We carry out the computational details in one
further case of interest (Theorem D).

\begin{proof}[The proof of Theorem D]
We have a basis
\[\mathcal B =\{t^i - t^{i+2q}\mid (i,2q)=1, i\equiv 1\Mod{4}, 1\leqq i < 4q\}\]
for $\rtildefree{G}$, so it remains to work out the relations
given by topological similarity.
Notice that $\Rtildefree_h(G) = \Rtildefree_n(G)$, and that
the sum of any two elements in $\mathcal B$ lies in $\Rtildefree_n(G)$.
Moreover, by Corollary \ref{sixdim} there are no $6${--}dimensional similarities
for $G$.

Now suppose that $V_1\oplus W\oplus \rminus \oplus \rplus
\sim_t V_2\oplus W\oplus \rminus \oplus \rplus$, for some
complex $G${--}representation $W$.
Then $\Res_H V_1 \cong \Res_H V_2$ since $q$ is odd, and
by Theorem C we get a similarity of the form
$V_1\oplus \Rminus{l} \oplus
 \rplus \sim_t V_2\oplus \Rminus{l} \oplus \rplus$.
But by Theorem \ref{sixteenone} this implies that 
$V_1\oplus \rminus\oplus \rplus
\sim_t V_1\oplus \rminus\oplus \rplus$.
Therefore, for any element $a_i=(t^i -t^{i+2q}) \in \mathcal B$
we have $2a_i \notin \Rfree_t(G)$ but $4a_i \in \Rfree_t(G)$.
However, in Section \ref{sixdiml} we determined the bounded
surgery obstructions for all these elements.
Since $i_{\ast}(u_{1,j}) = \langle -1\rangle \in \Coh{1}{\Units{A}}$
for all $j$ with $(j,2q)=1$
by Lemma \ref{seventeenone}, there are further stable relations
$a_1 + a_j \approx a_1 + a_k$, or $a_1 \approx a_j$
for all $j$.
It follows that a basis for $\rtildetopfree{G}$ is
given by $\{ a_1\mid 4a_1 \approx 0\}$.
\end{proof}
\newpage

\begin{thebibliography}{10}

\bibitem{cs1}
S.~E. Cappell and J.~L. Shaneson, \emph{Non-linear similarity}, Ann. of Math.
  (2) \textbf{113} (1981), 315--355.

\bibitem{cs2}
\bysame, \emph{The topological rationality of linear representations}, Inst.
  Hautes {\'E}tudes Sci. Publ. Math. \textbf{56} (1983), 309--336.

\bibitem{cs5}
\bysame, \emph{Torsion in {$L$}-groups}, Algebraic and Geometric Topology,
  Rutgers 1983, Lecture Notes in Mathematics, vol. 1126, Springer, 1985,
  pp.~22--50.

\bibitem{cssww1}
S.~E. Cappell, J.~L. Shaneson, M.~Steinberger, S.~Weinberger, and J.~West,
  \emph{The classification of non-linear similarities over {$\mathbb Z/2^r$}},
  Bull. Amer. Math. Soc. (N.S.) \textbf{22} (1990), 51--57.

\bibitem{cssw1}
S.~E. Cappell, J.~L. Shaneson, M.~Steinberger, and J.~West, \emph{Non-linear
  similarity begins in dimension six}, J. Amer. Math. Soc. \textbf{111} (1989),
  717--752.

\bibitem{cap1}
M.~C\'ardenas and E.~K. Pedersen, \emph{On the {K}aroubi filtration of a
  category}, K-theory \textbf{12} (1997), 165--191.

\bibitem{cp1}
G.~Carlsson and E.~K. Pedersen, \emph{Controlled algebra and the {N}ovikov
  conjectures for {K}- and {L}-theory}, Topology \textbf{34} (1995), 731--758.

\bibitem{carter1}
D.~Carter, \emph{Localizations in lower algebraic {K}-theory}, Comm. Algebra
  \textbf{8} (1980), 603--622.

\bibitem{dre1}
A.~W.~M. Dress, \emph{Induction and structure theorems for orthogonal
  representations of finite groups}, Ann. of Math. (2) \textbf{102} (1975),
  291--325.

\bibitem{hp1}
I.~Hambleton and E.~K. Pedersen, \emph{Bounded surgery and dihedral group
  actions on spheres}, J. Amer. Math. Soc. \textbf{4} (1991), 105--126.

\bibitem{hp2}
\bysame, \emph{Non-linear similarity revisited}, Prospects in Topology:
  (Princeton, NJ, 1994), Annals of Mathematics Studies, vol. 138, Princeton
  Univ. Press, Princeton, NJ, 1995, pp.~157--174.
  
\bibitem{hp3}  
\bysame, \emph{Topological equivalence of linear representations
                for cyclic groups}: \partone, Ann. of Math. (to appear).


\bibitem{htw3}
I.~Hambleton, L.~Taylor, and B.~Williams, \emph{An introduction to the maps
  between surgery obstruction groups}, Algebraic Topology (Aarhus, 1982),
  Lecture Notes in Mathematics, vol. 1051, Springer, Berlin, 1984, pp.~49--127.

\bibitem{ht1}
I.~Hambleton and L.~R. Taylor, \emph{A guide to the calculation of surgery
  obstruction groups}, Surveys in Surgery Theory, Volume I, Annals of
  Mathematics Studies, vol. 145, Princeton Univ. Press, 2000, pp.~225--274.

\bibitem{hsp1}
W-C. Hsiang and W.~Pardon, \emph{When are topologically equivalent
  representations linearly equivalent}, Invent. Math. \textbf{68} (1982),
  275--316.

\bibitem{mr1}
I.~Madsen and M.~Rothenberg, \emph{On the classification of {$G$}-spheres {I}:
  equivariant transversality}, Acta Math. \textbf{160} (1988), 65--104.

\bibitem{milgran2}
R.~J. Milgram and A.~A. Ranicki, \emph{The {$L$}-theory of {L}aurent extensions
  and genus {$0$} function fields}, J. Reine Angew. Math. \textbf{406} (1990),
  121--166.

\bibitem{o1}
R.~Oliver, \emph{Whitehead {G}roups of {F}inite {G}roups}, London Math. Soc.
  Lecture Notes, vol. 132, Cambridge Univ. Press, 1988.

\bibitem{pw1}
E.~K. Pedersen and C.~Weibel, \emph{A nonconnective delooping of algebraic
  {K}-theory}, Algebraic and Geometric Topology, (Rutgers, 1983), Lecture Notes
  in Mathematics, vol. 1126, Springer, Berlin, 1985, pp.~166--181.

\bibitem{ra13}
A.~A. Ranicki, \emph{The double coboundary}, Letter to B. Williams, available
  at the web page http://www.maths.ed.ac.uk/~aar/surgery/brucelet.pdf, 1981.

\bibitem{ra8}
\bysame, \emph{Exact {S}equences in the {A}lgebraic {T}heory of {S}urgery},
  Math. Notes, vol.~26, Princeton Univ. Press, 1981.

\bibitem{ra5}
\bysame, \emph{Additive {$L$}-theory}, K-theory \textbf{3} (1989), 163--195.

\bibitem{ra7}
\bysame, \emph{Algebraic {L}-theory and {T}opological {M}anifolds}, Cambridge
  Tracts in Math., vol. 102, Cambridge Univ. Press, 1992.

\bibitem{drh2}
G.~de Rham, \emph{Sur les nouveaux invariants topologiques de {M}
  {R}eidemeister}, Mat. Sbornik \textbf{43} (1936), 737--743, Proc.
  International Conference of Topology (Moscow, 1935).

\bibitem{drh1}
\bysame, \emph{Reidemeister's torsion invariant and rotations of {$S^n$}},
  Differential Analysis, (Bombay Colloq.), Oxford University Press, London,
  1964, pp.~27--36.

\bibitem{se1}
J.-P. Serre, \emph{Linear {R}epresentations of finite groups}, Graduate Texts
  in Mathematics, vol.~42, Springer, 1977.

\bibitem{sw1}
R.~Swan, \emph{Induced representations and projective modules}, Ann. of Math.
  (2) \textbf{71} (1960), 267--291.

\bibitem{swe1}
R.~G. Swan and E.~G. Evans, \emph{K-theory of finite groups and orders},
  Lecture Notes in Mathematics, vol. 149, Springer, 1970.

\bibitem{w2}
C.~T.~C. Wall, \emph{Surgery on {C}ompact {M}anifolds}, Academic Press, New
  York, 1970.

\bibitem{w8}
\bysame, \emph{On the classification of {H}ermitian forms. {II}. {S}emisimple
  rings}, Invent. Math. \textbf{18} (1972), 119--141.

\bibitem{yo1}
C.~Young, \emph{Normal invariants of lens spaces}, Canad. Math. Bull.
  \textbf{41} (1998), 374--384.

\end{thebibliography}

\providecommand{\bysame}{\leavevmode\hbox to3em{\hrulefill}\thinspace}
\providecommand{\MR}{\relax\ifhmode\unskip\space\fi MR }
\providecommand{\MRhref}[2]{%
  \href{http://www.ams.org/mathscinet-getitem?mr=#1}{#2}
}
\providecommand{\href}[2]{#2}

\end{document}